\documentclass[11pt]{article}
\def\version{4.10.2018}\def\users{}  %
\newcommand\COL[1]{{\color{black}#1}} 
\newcommand\REF[1]{{\color{black}#1}} 
\newcommand\REFF[1]{{\color{black}#1}} 
\def\users{final-layout}  
\usepackage{times}
\usepackage{enumerate}
\usepackage{amsmath,amsfonts,amssymb,color}
\usepackage{mathrsfs} 
\usepackage{epsfig}
\usepackage{psfrag}
\usepackage{upgreek} 
\usepackage{cite}

\newtheorem{theorem}{Theorem}[section]

\newtheorem{proposition}[theorem]{Proposition}

\newtheorem{remark}[theorem]{Remark}

\numberwithin{equation}{section}

\usepackage{ifthen}
\ifthenelse{\equal{\users}{final-layout}}{
\pagestyle{myheadings}
\markboth{T.Roub\'\i\v cek
}{A general model for various adhesive frictional contacts 
at small strains}
}{
\usepackage{fancyhdr}
\pagestyle{fancy}
\headheight=28pt\headwidth=16.5cm
\definecolor{gray}{gray}{0.5}
\rhead{\color{gray}Adhesive frictional contacts\\
T.Roub\'\i\v cek, 
}
\chead{}
\lhead{Version \version, file: \jobname.tex \\
compiled:
\number\day.\number\month.\number\year\ at
\the\hour:\ifnum\minute<10 0\fi\the\minute\ h }
}

\ifthenelse{\equal{\users}{final-layout}}{}{
\usepackage[notcite,notref,color]{showkeys}
\definecolor{labelkey}{rgb}{1.,.2,0.}

}

\newcount\hour \newcount\minute
\hour=\time
\divide \hour by 60
\minute=\time
\loop \ifnum \minute > 59 \advance \minute by -60 \repeat

\usepackage[normalem]{ulem}
\definecolor{brown}{rgb}{0.5,0,0}
\ifthenelse{\equal{\users}{final-layout}}{

    \newcommand{\DELETE}[1]{}
    
    \newcommand{\COMMENT}[1]{}
    \newcommand{\COLOR}[2][black]{{\color{#1}{#2}}}
    \newcommand{\TINY}[1]{}
    \newcommand{\MARGINOTE}[1]{}
}
{

 \newcommand{\DELETE}[1]{{\color{brown}\sout{#1}\color{black}}}
 
 \newcommand{\COMMENT}[1]{{\color{red}\uuline{#1}\color{black}}}
 \newcommand{\COLOR}[2][black]{{\color{#1}{#2}}}
 \newcommand{\TINY}[1]{{\tiny#1}}
 \newcommand{\MARGINOTE}[1]{\marginpar{\color{red}\tiny\texttt{#1}}}
}

\newcommand{\Item}[3][0.06]{\parbox[t]{#1\textwidth}{#2}\hfill%
      \parbox[t]{\dimexpr\textwidth-#1\textwidth}{#3}\vspace*{.8mm}}
\renewcommand\dot[1]{\mathchoice
                 {{\buildrel{\hspace*{.1em}\text{\LARGE.}}\over{#1}}}
                 {{\buildrel{\hspace*{.1em}\text{\Large.}}\over{#1}}}
                 {{\buildrel{\hspace*{.1em}\text{\large.}}\over{#1}}}
                 {{\buildrel{\hspace*{.1em}\text{\large.}}\over{#1}}}}
\newcommand\DT{\dot}
\newcommand\DDT[1]{\mathchoice
   {{\buildrel{\hspace*{.1em}\text{\LARGE.\hspace*{-.1em}.}}\over{#1}}}
   {{\buildrel{\hspace*{.1em}\text{\Large.\hspace*{-.1em}.}}\over{#1}}}
   {{\buildrel{\hspace*{.1em}\text{\large.\hspace*{-.1em}.}}\over{#1}}}
   {{\buildrel{\hspace*{.1em}\text{\large.\hspace*{-.1em}.}}\over{#1}}}}
\newcommand\JUMP[3][]{\mathchoice
                  {\big[\hspace*{-.3em}\big[#2\big]\hspace*{-.3em}\big]_{#3}^{#1}}
                   {[\hspace*{-.15em}[#2]\hspace*{-.15em}]_{#3}^{#1}}
                   {[\![#2]\!]_{#3}^{#1}}
                   {[\![#2]\!]_{#3}^{#1}}}
\newcommand{\wt}[1]{\mathchoice
     {\text{\small$\widetilde{\text{\normalsize$#1$}}\hspace*{.03em}$}}
                    {\text{\small$\widetilde{\text{\normalsize$#1$}}$}}
                    {\widetilde{#1\hspace*{-.02em}}\hspace*{.03em}}
                    {\widetilde{#1}}}

\def\R{{\mathbb R}}
\newcommand\bbC{\mathbb C}
\newcommand\bbD{\mathbb D}
\newcommand\bbE{\mathbb E}
\newcommand\bbK{\mathbb K}
\newcommand\bbM{\mathbb M}
\newcommand\bbB{\mathbb B}
\newcommand\SYM{\R_\mathrm{sym}^{d\times d}}
\newcommand\GC{\mathchoice{\varGamma_{\hspace*{-.15em}\mbox{\tiny\rm C}}}
                          {\varGamma_{\hspace*{-.15em}\mbox{\tiny\rm C}}}
                          {\varGamma_{\hspace*{-.05em}\mbox{\tiny\rm C}}}
                          {\varGamma_{\hspace*{-.05em}\mbox{\tiny\rm C}}}}
\newcommand{\Sdir}{\Sigma_{\mbox{\tiny\rm D}}}
\newcommand{\Gdir}{\varGamma_{\mbox{\tiny\rm D}}}
\newcommand{\udir}{u_{\mbox{\tiny\rm D}}}
\newcommand{\barudir}{\bar u_{\mbox{\tiny\rm D}}}
\newcommand{\DTbarudir}{\DT{\bar u}_{\mbox{\tiny\rm D}}}
\newcommand{\DDTbarudir}{\DDT{\bar u}_{\mbox{\tiny\rm D}}}
\newcommand{\Snew}{\Sigma_{\mbox{\tiny\rm N}}}
\newcommand{\Gnew}{\varGamma_{\mbox{\tiny\rm N}}}
\newcommand\SC{\Sigma_{\mbox{\tiny\rm C}}}

\newcommand{\nablaS}{\nabla_{\scriptscriptstyle\textrm{\hspace*{-.3em}S}}^{}}
\newcommand{\divS}{\mathrm{div}_{\scriptscriptstyle\textrm{\hspace*{-.1em}S}}^{}}
\newcommand{\thetaS}{\theta_{\scriptscriptstyle\textrm{\hspace*{-.1em}A}}^{}}
\newcommand{\thetaStauk}{\theta_{\scriptscriptstyle\textrm{\hspace*{-.1em}A},\tau}^{k}}
\newcommand{\underthetaStau}{\underline\theta_{\scriptscriptstyle\textrm{\hspace*{-.1em}A},\tau}^{}}
\newcommand{\underthetaBtau}{\underline\theta_{\scriptscriptstyle\textrm{\hspace*{-.1em}B},\tau}^{}}
\newcommand{\thetaStaukk}{\theta_{\scriptscriptstyle\textrm{\hspace*{-.1em}A},\tau}^{k-1}}
\newcommand{\thetaBtauk}{\theta_{\scriptscriptstyle\textrm{\hspace*{-.1em}B},\tau}^{k}}
\newcommand{\thetaBtaukk}{\theta_{\scriptscriptstyle\textrm{\hspace*{-.1em}B},\tau}^{k-1}}
\newcommand{\thetaStaull}{\theta_{\scriptscriptstyle\textrm{\hspace*{-.1em}A},\tau}^{l-1}}

\newcommand{\gammaC}{\gamma_{\scriptscriptstyle\textrm{\hspace*{-.1em}C}}^{}}
\newcommand{\gammaCprime}{\gamma_{\scriptscriptstyle\textrm{\hspace*{-.1em}C}}'}
\newcommand{\etaS}{\eta_{\scriptscriptstyle\textrm{\hspace*{-.05em}A}}^{}}
\newcommand{\indexS}{_{\scriptscriptstyle\textrm{\hspace*{.05em}A\hspace*{-.2em}}}^{}}
\newcommand{\indexV}{_{\scriptscriptstyle\textrm{\hspace*{.05em}B\hspace*{-.1em}}}^{}}
\newcommand{\zetaS}{\zeta_{\scriptscriptstyle\textrm{\hspace*{-.1em}A}}^{}}
\newcommand{\zetaSeq}{\zeta_{\scriptscriptstyle\textrm{\hspace*{-.1em}A}}^{\scriptscriptstyle\textrm{\hspace*{-.1em}eq}}}
\newcommand{\zetaBeq}{\zeta_{\scriptscriptstyle\textrm{\hspace*{-.1em}B}}^{\scriptscriptstyle\textrm{\hspace*{-.1em}eq}}}
\newcommand{\DTzetaS}{\DT\zeta_{\scriptscriptstyle\textrm{\hspace*{-.1em}A}}^{}}
\newcommand{\muS}{\mu_{\scriptscriptstyle\textrm{\hspace*{-.1em}A}}^{}}
\newcommand{\KS}{K_{\scriptscriptstyle\textrm{\hspace*{-.1em}A}}^{}}
\newcommand{\MS}{M_{\scriptscriptstyle\textrm{\hspace*{-.1em}A}}^{}}
\newcommand{\betaS}{\beta_{\scriptscriptstyle\textrm{\hspace*{-.1em}A}}^{}}
\newcommand{\thetaB}{\theta_{\scriptscriptstyle\textrm{\hspace*{-.1em}B}}^{}}
\newcommand{\zetaB}{\zeta_{\scriptscriptstyle\textrm{\hspace*{-.1em}B}}^{}}
\newcommand{\etaB}{\eta_{\scriptscriptstyle\textrm{\hspace*{-.05em}B}}^{}}
\newcommand{\muB}{\mu_{\scriptscriptstyle\textrm{\hspace*{-.1em}B}}^{}}

\newcommand{\MB}{M_{\scriptscriptstyle\textrm{\hspace*{-.1em}B}}^{}}
\newcommand{\betaB}{\beta_{\scriptscriptstyle\textrm{\hspace*{-.1em}B}}^{}}
\newcommand{\bbMB}{{\mathbb M}_{\scriptscriptstyle\textrm{\hspace*{-.1em}B}}^{}}
\newcommand{\bbMS}{{\mathbb M}_{\scriptscriptstyle\textrm{\hspace*{-.1em}A}}^{}}
\renewcommand\d{\mathrm d}
\newcommand{\one}{\mathbf{1}}
\newcommand{\K}{\mathbf{k}}
\newcommand{\bftheta}{\boldsymbol\theta}
\newcommand{\bfzeta}{\boldsymbol\zeta}
\newcommand{\bfeta}{\boldsymbol\eta}
\newcommand{\bfpsi}{\boldsymbol\psi}
\newcommand{\bfvartheta}{\boldsymbol\vartheta}
\newcommand{\bfmu}{\boldsymbol\mu}

\newcommand{\bulet}{{\hspace{.1em}\scriptstyle\bullet\hspace{.1em}}}

\allowdisplaybreaks 

\begin{document}
\def\rmn{{_{\rm N}}}
\def\rmt{{_{\rm T}}}


\newcommand{\subjclass}[1]{\bigskip\noindent\emph{2010 Mathematics Subject Classification:}\enspace#1}
\newcommand{\keywords}[1]{\noindent\emph{Keywords:}\enspace#1}




\title{\bf
\REF{A general} \COL{thermodynamical} model for adhesive\\frictional contacts 
between
\COL{viscoelastic or}\\ \COL{poro-viscoelastic} bodies at small strains}

\author{{\sc Tom\'{a}\v{s} Roub\'\i\v{c}ek}\\[1em]
Institute of Thermomechanics, Czech Academy of Sciences\\Dolej\v skova~5,
CZ-182~08 Praha 8, Czech Republic\\
tomas.roubicek@mff.cuni.cz
}

\date{}

\maketitle

\begin{abstract}
A general model covering a large variety of the so-called adhesive 
or cohesive\COL{, possibly also frictional,} contact \COLOR{interfaces} between 
visco-elastic bodies 
 with 
\COL{inertia} considered in a 
thermodynamical context is presented. A semi-implicit time discretisation
which \COL{conserves energy,} is numerically stable and convergent\COL{,} 
and which advantageously decouples
the system is devised.
\REF{An extension to porous media with adhesive contact influenced by 
a diffusant concentration is devised, too.}



\medskip

\end{abstract}

\def\friction{\mathfrak{f}}

\section{Introduction}

Contact mechanics is an important and widely developed part of continuum
mechanics of solids, intensively studied during past decades from the
modelling, computational, and analytical aspects, too~\cite{SoHaSh06AACP,SofMat12MMCM,SoMa09MMCM,Wrig06CCM,Ya13NMCM}. We focus on
a class of models for contacts with adhesion (cohesion) and friction in the
thermodynamical context. Our goal is to formulate a rather general
multi-purpose model unifying a relatively 
\COL{wide} \REFF{menagerie} of such particular
adhesion-friction models and to discuss its thermodynamics, outline
its \COL{conceptually efficient} discretisation \COL{amenable for 
calculation of vibrations/waves}.

\COLOR{
It is certainly impossible to present in \REFF{a} comprehensive way the relevant
literature to the part of contact mechanics covered by the general model
developed in the present paper. Anyhow, from the engineering point of view,
a detailed  analysis of related variational formulations of 
interface
cohesive models  was developed in~\cite{DeP14AVAT}.  Another  related and 
thermodynamically based approach for general imperfect interfaces was proposed 
and numerically studied in~\cite{EsJaSt16IJSS,JaMBSt13TSLD,JaKaSt14CMAME}, 
considering e.g.\ an interface temperature different from the temperatures of 
adjacent bulks, \COL{cf.\ also \cite{BoBoRo09TEAC,BoBoRo11LTBT,BoBoRo15MEBA} 
for analytical considerations}. Nevertheless, \COL{except \cite{BoBoRo15MEBA},}
friction dissipation was not 
taken into account in these works. Coupling of cohesive models and friction 
contact was developed after the pioneering  work~\cite{Tv90MSE}, which assumed  
the friction is activated after complete cohesive failure, in a series of 
relevant contributions~\cite{AlSa06IJNME,BiMr05IJSS,ChGiLe97IJDM,DelPiRa10EJM,GuiNgu14IJSS, LiGeSo01IJSS,LiGuPe16IJF,PaMaBo2016EFM,RaCaCo99CMCA,SeSaAl15EJM,SnMo13IJNME,TaCu03EJM,vdMeSl13C} among others, providing, for isothermal 
configurations, a continuous transition from the cohesive state to the frictional contact state.
However, no rigorous proof of the solution existence or the convergence of a numerical scheme to the exact solution was provided in these works.
From the mathematical point of view, overviews of existing partial models in 
literature can be found, e.g., without Coulomb friction ($\friction=0$) 
in \cite{ChShSo04DFCA,RoSoVo13MBMF,ScaSch17CPVB} (isothermal), without
inertia/friction and
interfacial plasticity but with fracture-mode-mixity sensitive \REF{debonding 
(sometimes called also} delamination) in \cite{KrPaRo15QACD} (isothermal) 
or \cite{RosRou11TARI,RosRou13ACDM} (anisothermal); see also a review chapter \cite{RoKrZe14Chapter}.}
%
%
%
%
%
%
%
%
%
%
%


Models of adhesive contacts are to a great extent inspired
by bulk damage possibly combined with plasticity. The role of a scalar
damage variable is then replaced by a\COLOR{n
interface} damage (\REF{a concept invented by M.\,Fr\'emond 
\cite{Frem85DAS,Frem02NST}}, also
called \REF{debonding} or delamination\COLOR{, the term adopted hereinafter}) 
variable and plasticity can be transferred to a plastic slip, 
\REF{cf.\ Remark~\ref{rem-pi}
for more details.} Like in the bulk models,
there are many options:\\
\Item{(i)}{delamination can influence the elasticity
of the adhesive (through decaying elastic moduli) or
the plasticity (through decaying plastic yield stress for the plastic slip),
or both;}
\Item{(ii)}{vice versa, plasticity can influence delamination
indirectly through influencing the stress and displacement jump,
or directly through influencing activation threshold for delamination;}
\Item{(iii)}{delamination evolution can be either unidirectional
or reversible (i.e.\ admitting a so-called healing);}
\Item{(iv)}{delamination and interfacial plasticity can be considered
rate-independent
(and then choosing an appropriate concept of solution becomes a vital part
of the model itself and should carefully be considered, cf.\ \cite{MieRou15RIS})
or rate-dependent (visco-plasticity, viscous delamination),
i.e.\ 4 options altogether, or more if also healing is
involved either as rate-dependent or rate-independent;}
\Item{(v)}{interfacial plasticity can be either with hardening
or without hardening;}
\Item{(vi)}{length scale (gradients) can be considered in \REF{interfacial}
damage or/and plasticity.}
In contrast to the mentioned bulk models, the interfacial analog
can additionally be combined with friction. We can then identify
other options (beside the trivial option just ignoring friction):\\
\Item{(vii)}{friction acts on the displacement velocity, or can be realized
through the interfacial plastic slip;}
\Item{(viii)}{the contact is unilateral or bilateral, in the former case
the adhesion being mode sensitive or insensitive, possibly in
the normal-compliance variant counting with a small penetration
(microscopically related to various
asperities, cf.\ e.g.\ \cite{TCOY98CMMM} for a detailed discussion).}
Of course, as a general feature in (contact) mechanics, there is still
an option as far as temperature variation and related heat transfer concerns:\\
\Item{(ix)}{either to ignore it by considering the model isothermal, or
to consider heat transfer in the bulk and through the contact 
interface
but ignoring the heat capacity and the heat transfer in the adhesive
itself, or even to involve this latter heat problem in the adhesive, too.}
In addition, various non-linearities can be considered at various
parts of the models, leading e.g.\ to various cohesive-contact models etc.
Eventually, there are standard options as far as the bulk models:
purely elastic or visco-elastic with various rheologies (possibly involving
creep) or some additional inelastic processes like plasticity,
purely quasistatic models versus inertial effects involved,
simple or nonsimple materials, cf.\ \cite{RosRou13ACDM}, linear or nonlinear
materials, purely isothermal versus anisothermal with or without thermal
expansion, etc.

Contact mechanics has many applications in engineering or in (for example geo-)
physics, and various combinations of the above identified options (i)--(ix)
may serve for various purposes. Of course, not all combinations give sensible
models or are amenable for a rigorous analysis supporting
also executable numerical schemes with guaranteed stability and convergence.

The above outlined rich menagerie of inelastic-contact models however should not
create a false impression that it covers everything. In particular, in this
paper we ignore phenomena like fatigue or wear, and \COL{we 
confine ourselves} to small strains and, in addition, we also assume small 
displacements (and thus a prescribed contact zone, avoiding e.g.\ contact 
of rolling bodies \COL{or a varying Hertz' contact zone}).

\COL{Moreover, there is also a lot of contact-mechanical models based on
hemivariational inequalities involving nonsmooth nonconvex functions, 
cf.\ e.g.\ \cite{Bart09HIMD,BaHaMo05MMDM,MiOcSo13NIHI}. Often, this
mathematically involved and numerically more difficult setting results 
rather from lacking of suitable internal variables. In contrast, the 
definite advantage of our general model is that it relies on classical 
variational inequalities involving convex functions only.}

The plan of this paper is the following: In Section~\ref{sect-model}
we formulate the general model \COL{coupling the adhesive layer with 
the viscoelastic bulk. Then,} in Section~\ref{sec-thermodyn}
we \COL{devise its ``monolithic'' formulation treating the bulk and
the adhesive in a unified way, and} discuss its energetics and thermodynamical
consistency. Then, in Section~\ref{sect-disc}, a semi-implicit 
\COL{energy-conserving} time
discretisation allowing for efficient computer implementation is devised
and its numerical stability and convergence is outlined. \REF{Eventually, 
in Section~\ref{poro}, we outline the extension towards poro-viscoelastic
bulk materials and adhesives, devising again a monolithic formulation.}
For readers' convenience, let us summarize
the basic notation used in what follows in Table~\ref{Tab_Notation}.

\begin{table}[ht]
\centering
\fbox{
\begin{minipage}[t]{17em}
\small

$\varOmega,\varOmega_1,\varOmega_2\subset\R^d$ domains, 
$\ \varOmega_1{\cap}\varOmega_2=\emptyset$,

$d=2,3$ dimension of the problem,

\COLOR{$\GC$ contact interface, 
$\ \varOmega\!\setminus\!\GC=\varOmega_1{\cup}\varOmega_2$}

$\SYM:=\{A\!\in\!\R^{d\times d};\ A=A^\top\}$,


$u:Q{\setminus}\SC\to\R^d$ displacement,

$\pi{:}\SC\to\R^{d-1}$ plastic-like interfacial slip,


$\alpha:\SC\to[0,1]$ interfacial damage, 

\hspace*{6.5em}(delamination) variable,

$\thetaS\!:\SC\to
\R^+$
temperature \COL{in the adhesive},

$\thetaB:
Q{\setminus}\SC{\to}
\R^+$
bulk temperature,

$e(u)
=\frac12\nabla u^\top\!+\frac12\nabla u$
%
small-strain tensor,

$\JUMP{u}{}:\SC\to\R^d$ jump of $u$ across $\GC$,

\COLOR{$\JUMP{u}{\rmn}$ normal jump of $u$ on $\GC$,}

\COLOR{$\JUMP{u}{\rmt}$ tangential jump of $u$ on $\GC$,}

\COLOR{$\boldsymbol{E}:u\mapsto(e(u),\JUMP{u}{})$,}

$\etaS\!:\SC\to
\R$ 
entropy \COL{in the adhesive},

$\etaB:Q{\setminus}\SC\to
\R$ specific bulk entropy,

\COL{$\vartheta\indexS$ heat part of internal energy in the adhesive,}

\COL{$\vartheta\indexV$ heat part of the bulk internal energy,}



$f:\Snew\to\R^d$  applied traction force,

$g:Q\to\R^d$  applied bulk force (gravity),



\REF{$\sigma$ stress in the bulk},

\REF{$\zetaS$ the content of diffusant in the adhesive},

\REF{$\zetaB$ the content of diffusant in the bulk},

$\bbC
\in\R^{d^4}$ elasticity tensor,

$\bbD
\in\R^{d^4}$ viscosity tensor,

$\bbK\indexS
\in\R^{\COL{(d{-}1)\times(d{-}1)}}$ 
adhesive heat conductivity,

$\bbK\indexV=\bbK\indexV(\thetaB)
\in\R^{\COL{d\times d}}$ bulk heat conductivity,

\end{minipage}\ \
\begin{minipage}[t]{19em}
\small

$\bbE
\in\R^{d\times d}$ thermal-expansion tensor,

$\varrho$ mass density,

\REF{$\bbMS$,\,$\bbMB$ mobility tensors in the adhesive and bulk}, 

\REF{$m$ mobility between the  adhesive and the bulk},

\REF{$\MS,\,\MB$ Biot moduli in the adhesive and bulk},

\REF{$\betaS,\,\betaB$ Biot coefficients in the  adhesive and bulk},

\REF{$\muS,\muB$ chemical potentials in the  adhesive and bulk},

\REF{$\sigma_{_{\rm C}}^{}$, $\sigma_{_{\rm F}}^{}$ stresses acting in dry-friction law,}

\REF{$\sigma_{_{\rm T}}^{},$ $\sigma_{_{\rm N}}^{}$ tangential/normal components of traction,}

$c\indexV=c\indexV(\thetaB)$ bulk heat capacity,

$c\indexS=c\indexS(\thetaS)$ adhesive heat capacity,

$\K\indexS(\REFF{\JUMP{u}{\rmn}},\alpha,\thetaS)$ bulk-to-adhesive
heat
 coefficients,




$\udir:\Sdir
\to\R^d$  prescribed
boundary displacement,

$\friction(\alpha,\thetaS)$ friction coefficient,



$\sigma_{\rm y}(\alpha,\thetaS)$ yield-stress for plastic slip,

$a_1(\REFF{\JUMP{u}{}},\alpha,\thetaS,\DT\alpha)$ delamination dissipation 
potential,

$a_0(\alpha)+b_0(\alpha)\thetaS$ \COLOR{stored} energy of delamination,

$\gamma_\rmt(\alpha,\cdot)$, $\gamma_\rmn(\alpha,\cdot)$
elastic moduli of the adhesive,

$\gammaC(\REFF{\JUMP{u}{}})$ normal-compliance potential,

$d_\rmn=d_\rmn(\alpha,\thetaS)$ normal viscosity of adhesive,

$d_\rmt=d_\rmt(\alpha,\thetaS)$ tangential viscosity of adhesive,

$\kappa_{\rm H}\ge0$ hardening coefficient for plastic slip,

$\kappa_1>0$ scale coefficient of the gradient of $\pi$,

$\kappa_2>0$ scale coefficient of the gradient of $\alpha$,

\COL{$\varPsi$ total free energy},

\COL{$\varXi$ total potential of dissipative forces},

\COL{$R$ total dissipation rate},
 
\COL{$M$ total kinetic energy}.

\end{minipage}\medskip
}
\caption{Summary of the basic notation used through the paper.
For the notation $Q$, $\SC$, $\Sdir$,\ see \eqref{notation} 
and Fig.~\ref{fig:geometry} below.
}\label{Tab_Notation}
\end{table}




\section{The general model}\label{sect-model}

\begin{figure}[ht]
\centering
\psfrag{O1}{\footnotesize $\varOmega_1$}
\psfrag{O2}{\footnotesize $\varOmega_2$}
\psfrag{GC}{\footnotesize $\GC$}
\psfrag{n}{\footnotesize $\vec{n}$}
\psfrag{nu}{\footnotesize $\nu\indexS$}
\psfrag{n1}{\footnotesize $\vec{n}_1$}
\psfrag{n2}{\footnotesize $\vec{n}_2$}
\psfrag{[u],pi,zeta,thetas}{\footnotesize
$\pi$, $\alpha$, $\thetaS$, \REF{$\zetaS$, $\muS$}}
\psfrag{u,theta}{\footnotesize $u$, $\thetaB$, \REF{$\zetaB$, $\muB$}}
\psfrag{bulk}{\footnotesize\begin{minipage}[c]{8em}\baselineskip=11pt visco-elastic\\bulk $\varOmega\!\setminus\!\GC$\\\hspace*{1em}$=\varOmega_1{\cup}\varOmega_2$\end{minipage}}
\psfrag{adhesive}{\footnotesize\hspace*{0em}
\begin{minipage}[c]{9.5em}\baselineskip=8pt a
frictional uni- \\-lateral con\-tact glued\\by a dama\-geable, heat- -con\-ductive,
infi\-ni\-te\-si\-mally thin adhesive \COLOR{layer} (=the interface 
$\GC$)\!\!\end{minipage}}
\includegraphics[width=150mm]{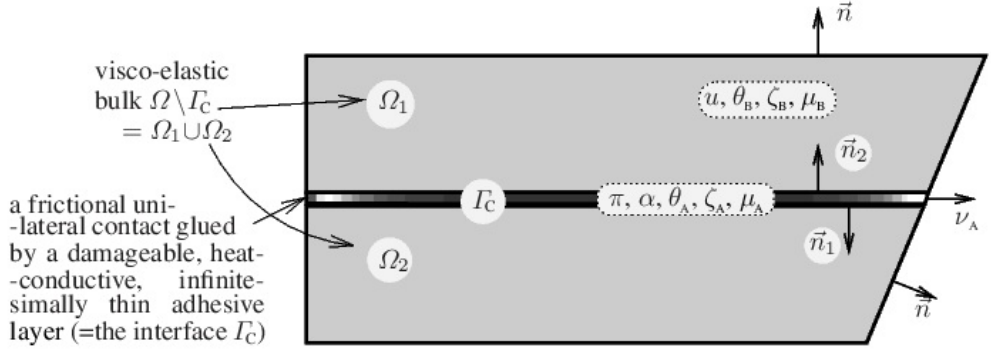}
\vspace*{-42em}
\vspace*{-.1em}
\caption{Schematic geometry of one contact interface $\GC$ across a 
2-di\-men\-sional domain $\varOmega$ 
and variables in the bulk and \COL{the adhesive} on $\GC$.
}\label{fig:geometry}
\end{figure}

We consider a bounded Lipschitz domain $\varOmega\subset\R^d$ split
to two (or more) bodies by some internal interface(s) $\GC$,
cf.\ Figure~\ref{fig:geometry}.
We use a concept of an ``interface'' temperature $\thetaS$ in an infinitesimally
thin but still heat conductive adhesive on the $(d{-}1)$-dimensional contact
\COLOR{interface} $\GC$ which glues two $d$-dimensional visco-elastic bodies,
cf.\ \cite{BoBoRo09TEAC,BoBoRo11LTBT,BouLan09NTDA} or
\cite[Sect.5.3.3.3]{MieRou15RIS} for such a model,
or in combination with a friction also \cite{BoBoRo15MEBA}.
For a relation to models which avoid this interface-temperature
concept see Remark~\ref{rem-no-thetaS}.
Rather for notational simplicity, we will consider no more complicated
rheology than the linear Kelvin-Voigt  one (so that, in particular, no
additional internal variables in the bulk will be considered).


\COLOR{The unit normal vector to $\varOmega_i$ $(i=1,2)$ is 
denoted by $\vec{n}_i$. Hence $\vec{n}_1=-\vec{n}_2$ on $\GC$, while we denote
just by $\vec{n}=\vec{n}_i$ on $\partial\varOmega_i{\setminus}\GC$, 
Fig.~\ref{fig:geometry}. A convention for definition of jumps across $\GC$ 
for traces of fields defined in $\varOmega_1$ and $\varOmega_2$ is used 
{as well as their normal and tangential components}, namely 
\begin{align}\label{convention}
\JUMP{u}{}\!:=u_1|_{\GC}^{}\!-u_2|_{\GC}^{},\ \ 
\JUMP{u}{\rmn}\!\!:=
\JUMP{u}{}{{\cdot}}\vec{n}_2=\big(u_2|_{\GC}^{}\!-u_1|_{\GC}^{}\big)
{{\cdot}}\vec{n}_1,
\ \ \JUMP{u}{\rmt}\!\!:=\JUMP{u}{}\!-\JUMP{u}{\rmn}\vec{n}_2.
\end{align}
Thus, e.g.\ 
$\JUMP{u}{\rmn}>0$ means that the boundaries of $\varOmega_1$ and $\varOmega_2$ at 
$\GC$ separate (we refer to such situation as opening) whereas for 
$\JUMP{u}{\rmn}<0$ overlapping of these subdomains takes place.}

We further consider the boundary of $\varOmega$ split to two parts,
$\Gdir$ and $\Gnew$, where the Dirichlet and the Neumann boundary
conditions will be prescribed. Considering  a fixed time horizon $T>0$,
we will use the notation
\begin{subequations}\label{notation}\begin{align}
&Q:=(0,T)  \times \varOmega, \qquad \Sigma : = (0,T)  \times \partial
\varOmega, \quad\\& \SC:= (0,T)  \times\GC , \quad\ \Sdir:= (0,T)
\times \Gdir , \quad\ \Snew:= (0,T) \times\Gnew.
\end{align}\end{subequations}
One of the main ingredient of the model will be
the force 
\REF{balance} on the contact interface $\GC$:
\begin{subequations}\label{system}
\begin{align}
\label{adhes-form-d1}
& \JUMP{\sigma}{}\vec{n}=0\ \ \ \ \ \text{ with }\
\sigma:=\bbD e(\DT{u})+\bbC (e(u){-}\bbE\thetaB\COLOR{{+}\bbE\theta_{_{\rm R}}})
&\text{on }\SC,\hspace{1.2em}
\\&\nonumber
|\sigma_{_{\rm F}}^{}|\le-\friction\sigma_{_{\rm C}}^{}\ \ \text{ with }\
\friction=\friction(\alpha,\thetaS),\ \ \
\sigma_{_{\rm C}}^{}:=\gammaCprime({\JUMP{u}{}}_\rmn),
\\&\nonumber\hspace{6em}\text{ and }\ \
\sigma_{_{\rm F}}^{}:=\sigma_\rmt^{}
-\big[\gamma_\rmt\big]_v'\big(\alpha,\JUMP{u}{\rmt}\!{-}\pi\big)
-d_\rmt(\alpha,\thetaS)\big(\JUMP{\DT u}{\rmt}\!{-}\DT\pi\big),
\\&\hspace{6em}\text{ and }\ \
\sigma_\rmt^{}:=
\sigma\COLOR{\vec{n}_2}-\sigma_\rmn^{}\COLOR{\vec{n}_2}\ \ \text{ with }\ \
\sigma_\rmn:=
\COLOR{\vec{n}_2^\top}\sigma
\COLOR{\vec{n}_2}
&\text{on }\SC,\hspace{1.2em}
\\[.2em]
&|\sigma_{_{\rm F}}^{}|<-\friction\sigma_{_{\rm C}}^{}\ \ \ \Rightarrow\ \ \ \JUMP{\DT u}{\rmt}=0
&\text{on }\SC,\hspace{1.2em}
\\[.2em]
&|\sigma_{_{\rm F}}^{}|=-\friction\sigma_{_{\rm C}}^{}\ \ \ \Rightarrow\ \ \ \exists\lambda\ge0:\ \
\sigma_{_{\rm F}}^{}=\lambda\JUMP{\DT u}{\rmt}
&\text{on }\SC,\hspace{1.2em}
\\&
\sigma_\rmn^{}-\big[\gamma_\rmn\big]_v'\big({\alpha,}\JUMP{u}{\rmn}\big)
-\gammaCprime\big({\JUMP{u}{}}_\rmn\big)
-d_\rmn(\alpha,\thetaS)\JUMP{\DT u}{\rmn}
=0&\text{on }\SC,\hspace{1.2em}
\intertext{\COLOR{where $\theta_{_{\rm R}}$ is a reference temperature 
(considered constant over $\varOmega$ for simplicity). 
F}urthermore, we consider the flow-rule for delamination as:}
\nonumber
&
%
\partial_{\DT\alpha}^{}a_1\big(\JUMP{u}{},\COLOR{\alpha,}\thetaS,\DT\alpha\big)
+[\gamma_\rmt]_\alpha'\big(\alpha,\JUMP{u}{\rmt}\!{-}\pi\big)
+[\gamma_\rmn]_\alpha'\big(\alpha,\JUMP{u}{\rmn}\big)
\\&\label{adhes-form-d6}
\qquad\qquad\qquad\qquad\qquad+N_{[0,1]}(\alpha)\ni\divS(\kappa_2\nablaS\alpha)
+a_0'(\alpha)+b_0'(\alpha)\thetaS
%
 &\text{on }\SC,\hspace{1.2em}
\intertext{with $N_{[0,1]}$ denoting the normal
cone to the interval $[0,1]$, i.e.\ the subdifferential in the sense of
convex analysis of the indicator function $\delta_{[0,1]}$ of $[0,1]$.
The other ingredient is the flow-rule for the interfacial plastic slip}
&\label{adhes-form-d7}
\DT\pi\in N_{|.|\le\sigma_{\rm y}(\alpha,\thetaS)}^{}\big(\divS(\kappa_1\nablaS\pi)
-[\gamma_{\rmt}]_v'(\alpha,\JUMP{u}{\rmt}\!{-}\pi)\big)
&\text{on }\SC,\hspace{1.2em}
\intertext{and moreover the interface heat-transfer equation}
&\nonumber
c\indexS(\thetaS)\DT{\thetaS\!}
-\divS(\bbK\indexS(\thetaS)\nablaS\thetaS)
=\K\indexS\big(\JUMP{u}{\rmn},\alpha,\thetaS\big)\cdot
\big(\thetaB|_{\GC}^{}\!{-}\thetaS\one\big)
+\friction(\alpha,\thetaS)\gammaCprime\big({\JUMP{u}{}}_\rmn\big)\big|\JUMP{\DT u}{\rmt}\big|\hspace*{-6em}
\\&\label{heat-GC}\hspace*{9em}
+\sigma_{\rm y}(\alpha,\thetaS)\big|\DT\pi\big|
+\big(\partial_{\DT\alpha}a_1(\JUMP{u}{},\COLOR{\alpha,}\thetaS,\DT\alpha)
-\thetaS b_0'(\alpha)\big)\DT\alpha &\text{on }\SC,\hspace{1.2em}
\end{align}
\end{subequations}
where $\thetaB|_{\GC}$ denotes the pair of the traces of a
possibly discontinuous temperature field $\thetaB$ from both sides of the
interface $\GC$, and $\K\indexS:=(k_1,k_2)$, and $\one:=(1,1)$.
Thus, the expression $\K\indexS\cdot(\thetaB|_{\GC}^{}\!{-}\thetaS\one)$
means 
{$k_{1}({\thetaB}_1|_{\GC}^{}{-}\thetaS)+
k_{2}({\thetaB}_2|_{\GC}^{}{-}\thetaS)$} with
${\thetaB}_1$ and ${\thetaB}_2$ denoting temperature
fields in the two domains adjacent to the interface $\GC$
and $k_1$ and $k_2$ the bulk-to-adhesive
heat
 coefficients from the corresponding bulk domains.

In (\ref{system}f-h), ``$\nablaS$'' denotes the gradient over
$(d{-}1)$-dimensional surface
\COLOR{(interface)} $\GC$, i.e.\ the tangential derivative
defined as 
{$\nablaS v=\nabla v-(\nabla v{\cdot}\vec{n})\vec{n}$} for $v$ defined 
\COLOR{around} $\GC$\COLOR{ with $\vec{n}=\vec{n}_1$ (or, equally, 
$\vec{n}=\vec{n}_2$); note that $\nabla v$ takes into account 
the profile of $v$ in a neighbourhood of $\GC$ although $\nablaS v$
eventually a.e.\ depends only on the trace on $\GC$}. 
Further the corresponding surface divergence ``$\divS$''
is the trace of $\nablaS$\COLOR{; thus $-\divS\nablaS$ is the so-called 
Laplace-Beltrami operator}.

In the bulk, we consider the standard thermo-visco-elasticity
with thermal expansion consisting 
\REFF{of} the
force balance coupled with the heat-transfer equation:

\begin{subequations}\label{bulk-model}
\begin{align}
\label{eq6:adhes-class-form1}& \varrho\DDT{u} -\mathrm{div}\,\sigma
=g\ \ \ \ \text{ with $\sigma$ as in \eqref{adhes-form-d1}}
&\text{in  }
Q{\setminus}\SC,
\\
\label{eq6:adhes-class-form1bis}
&c\indexV(\thetaB)\DT{\thetaB}-\mathrm{div}\big(\bbK\indexV(\thetaB)\nabla\thetaB\big)
=\bbD e(\DT{u}){:} e(\DT{u})-\thetaB\bbC\bbE{:}e(\DT{u})
&\text{in }Q{\setminus}\SC
\end{align}
\end{subequations}
with
$g:Q\to\R^d$ \COL{being} the applied bulk force (typically 
gravitational)\COLOR{;
let us say that we consider $\bbE$ constant over each subdomain $\varOmega_i$
so that the reference temperature $\theta_{_{\rm R}}$ occurs explicitly 
only in \eqref{system} but not in \eqref{eq6:adhes-class-form1}.}

Further, we supplement \eqref{bulk-model} with standard boundary conditions
counting also with the heat transfer from the adhesive:
\begin{subequations}\label{BC}
\begin{align}
\label{eq6:adhes-class-form2} &u=\udir(t) &\text{on }\Sdir,\hspace{1.2em}
\\\label{eq6:adhes-class-form3-bis}
&
\sigma\vec{n}=f(t)
\ \ \ \ \text{ with $\sigma$ again from \eqref{adhes-form-d1}}&\text{on }\Snew,\hspace{1.2em}
\\\label{eq6:adhes-class-form3}
&(\bbK\indexV(\thetaB)\nabla\thetaB)\vec{n}=h_{\rm ext}&\text{on }\Sigma,\hspace{1.6em}
\\\label{BC-heat-GC}
&(\bbK\indexV({\thetaB}_i)\nabla{\thetaB}_i)\vec{n}_i
=\K\indexS\big(\REFF{\JUMP{u}{\rmn}},\alpha,\thetaS\big)\big({\thetaB}_i|_{\GC}{-}\thetaS\big),
\ \ \ \ \COL{i=1,2,} &\text{on }\SC,\hspace{1.2em}
\intertext{where $f:\Snew\to\R^d$ is the applied traction,
$h_{\rm ext}:\Sigma\to\R$ is some external heat flux,
{and} $i$ in \eqref{BC-heat-GC} distinguishes the particular
subdomains adjacent to the interface $\GC$,
while \eqref{system} is supplemented with the boundary
conditions for $\pi$, $\alpha$, and $\thetaS$ on the $(d{-}2)$-dimensional
boundary of $\GC$:}
&\nablaS\pi{\cdot}\nu\indexS=0,\ \ \ \ \ \
\nablaS\alpha{\cdot}\nu\indexS=0,\ \ \ \text{ and }\ \ \
\nablaS\thetaS{\cdot}\nu\indexS=0
&\hspace*{-3em}\text{on }(0,T)\times\partial\GC\hspace*{2em}
\end{align}
\end{subequations}
with $\nu\indexS$ denoting the normal to the $(d{-}2)$-dimensional boundary
$\partial\GC$ of the interface $\GC$\COLOR{, cf.\ Figure~\ref{fig:geometry}}.

In addition, we consider an initial-value problem by prescribing the
following initial conditions
\begin{align}\label{IC}
u|_{t=0}^{}=u_0,\ \ \ \ \DT u|_{t=0}^{}=v_0,\ \ \ \ \pi|_{t=0}^{}=\pi_0,
\ \ \ \ \alpha|_{t=0}^{}=\alpha_0,
\ \ \ \ \thetaS|_{t=0}^{}=\theta_{{\scriptscriptstyle\textrm{\hspace*{-.05em}A}},0}^{}
,\ \ \ \ \thetaB|_{t=0}^{}={\thetaB}_{,0}.
\end{align}

\begin{figure}[ht]
\centering
\psfrag{gamma-t}{\footnotesize ${[\gamma_\rmt]'}_\REFF{\!\JUMP{u}{\rmt}}^{}\!\!(\alpha,\cdot)$}
\psfrag{gamma-n}{\footnotesize ${[\gamma_\rmn]'}_\REFF{\!\JUMP{u}{\rmn}}^{}\!\!(\alpha,\cdot)$}
\psfrag{u-t}{\footnotesize $\JUMP{u}{\rmt}$}
\psfrag{u-n}{\footnotesize $\JUMP{u}{\rmn}$}
\psfrag{u-n+t}{\footnotesize $\JUMP{u}{}$}
\psfrag{C}{\footnotesize $\bbC$}
\psfrag{D}{\footnotesize $\bbD$}
\psfrag{z}{\footnotesize $\pi$}
\psfrag{r}{\footnotesize $\varrho$}
\psfrag{e(u)}{\footnotesize $e(u)$}
\psfrag{normal compliance}{\footnotesize \begin{minipage}[c]{19em}\baselineskip=8pt normal-compliance\\[-.2em]approximation of 
Signoring contact\end{minipage}}
\psfrag{sigma-y}{\footnotesize $\sigma_{\rm y}(\alpha,\thetaS)$}
\psfrag{d-t}{\footnotesize $d_\rmt(\alpha,\thetaS)$}
\psfrag{d-n}{\footnotesize $d_\rmn(\alpha,\thetaS)$}
\psfrag{gamma-c}{\footnotesize $\gammaCprime$}
\psfrag{sigma-c}{\footnotesize
$\sigma_{_{\rm C}}^{}=\gammaCprime({\JUMP{u}{}}_\rmn)$}
\psfrag{sigma-n}{\footnotesize $\sigma_\rmn^{}\!=\vec{n}^\top\sigma\vec{n}$}
\psfrag{sigma-f}{\footnotesize $\sigma_{_{\rm F}}^{}$}
\psfrag{sigma-t}{\footnotesize $\sigma_\rmt^{}\!=\sigma\vec{n}-\sigma_\rmn^{}\vec{n}$}
\psfrag{mu}{\footnotesize $\friction(\alpha,\thetaS)$}
\psfrag{tangential}{\small \begin{minipage}[c]{6em}\baselineskip=8pt 
tangential\\[-.0em]direction\end{minipage}}
\psfrag{normal}{\small \begin{minipage}[c]{6em}\baselineskip=8pt 
normal\\[-.0em]direction\end{minipage}}
\psfrag{bulk2}{\small \begin{minipage}[c]{6em}\baselineskip=9pt visco-elastic bulk $\varOmega_2$\end{minipage}}
\psfrag{bulk1}{\small \begin{minipage}[c]{6em}\baselineskip=9pt visco-elastic bulk $\varOmega_1$\end{minipage}}
\psfrag{damageable adhes}{\footnotesize damageable visco-elastic adhesive}
\psfrag{damageable adhesive}{\footnotesize \begin{minipage}[c]{10em}\baselineskip=8pt damageable visco-\\[-.1em]\hspace*{.3em}-elastic adhesive\end{minipage}}
\psfrag{damageable plasticity}{\footnotesize \begin{minipage}[c]{10em}\baselineskip=8pt  interfacial plasticity\\[-.1em]\hspace*{.9em}(damageable)\end{minipage}}
\psfrag{k-H}{\footnotesize $\kappa_{\rm H}$}
\psfrag{hardening}{\footnotesize hardening}
\psfrag{Coulomb friction}{\footnotesize Coulomb\hspace{1.6em}friction}
\includegraphics[width=190mm]{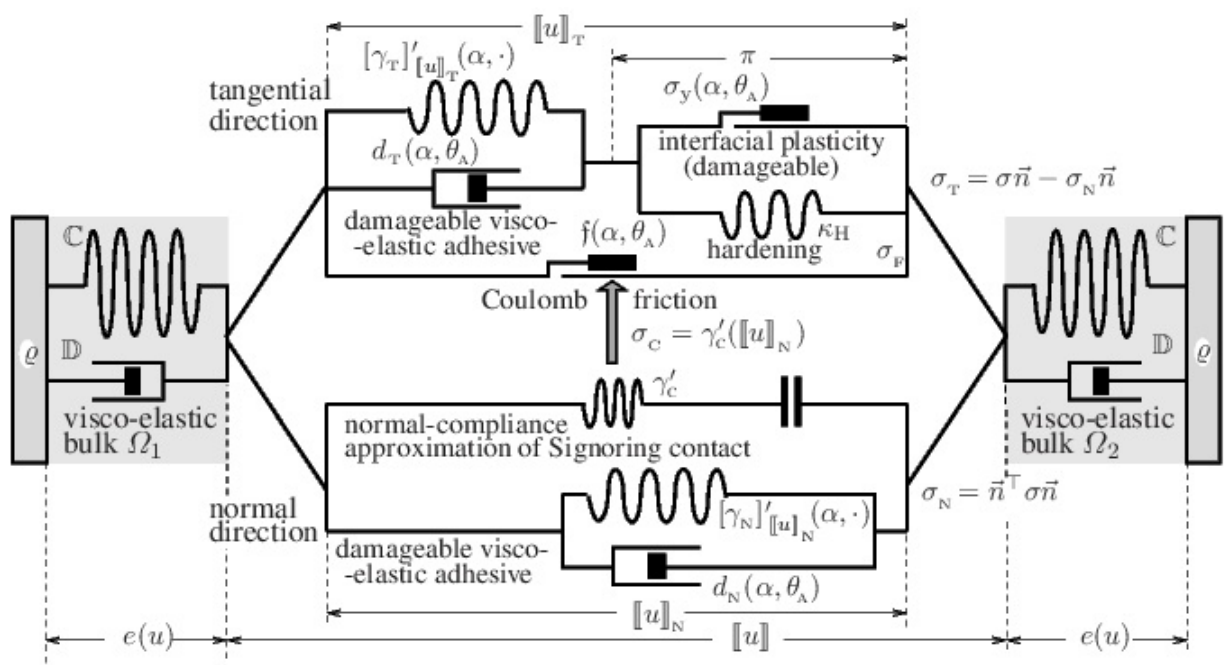}
\vspace*{-45em}
\vspace*{-.1em}
\caption{Schematic rheological model used on the \COL{adhesive} interface 
$\GC$ with one internal variable $\pi$ relevant for \COL{the} tangential 
\COL{slip} \REFF{and the unilateral contact for the normal displacement}. 
Thermal expansion in the bulk 
is not depicted, neither is evolution for delamination  $\alpha$ as the other
internal variable\COL{s} on $\GC$.
}\label{fig:rheology}
\end{figure}

The conceptual rheological diagram corresponding to the part of the model,
namely
\eqref{system}--\eqref{eq6:adhes-class-form1},
is shown in Fig.\,\ref{fig:rheology}.
The modelling assumption is that, under compression, the normal compliance 
dominates the response of the adhesive, i.e.\ $\gammaC(\REFF{\JUMP{u}{\rmn}}
)\gg\gamma_\rmn(\alpha,\REFF{\JUMP{u}{\rmn}})$ for
$\REFF{\JUMP{u}{\rmn}}<0$\COL{, otherwise rather $\sigma_{\rmn}^{}$ should enter
(\ref{system}b-d), cf.\ also Remark~\ref{rem-others}}.
\REFF{The standard symbols are used for elastic and viscous elements 
(springs and dashpots), respectively, and for the dry-friction-type
elements.}

\begin{remark}[More general $c\indexS$ and $\bbK\indexS$.]
\upshape
Some applications may need rather
$c\indexS=c\indexS(\alpha,\pi,\thetaS)$ and
$\bbK\indexS=\bbK\indexS(\alpha,\pi,\thetaS)$. The $(\alpha,\pi)$-dependent
interface heat capacity needs generalization of \eqref{def-of-C} \COLOR{below}
so that $\vartheta\indexS=C\indexS(\thetaS,\alpha,\pi)$ and yields some
adiabatic-like heat source/sink contribution to the heat equation, namely
$[C\indexS]_\alpha'(\thetaS,\alpha,\pi)\DT\alpha+
[C\indexS]_\pi'(\thetaS,\alpha,\pi)\DT\pi$.
\end{remark}

\begin{remark}[Models without interfacial temperature:
$c\indexS=0$ and $\bbK\indexS=0$.]\label{rem-no-thetaS}
\upshape
A particular case
arises when the adhesive has zero heat conductivity
$\bbK\indexS$ and capacity $c\indexS$. Then,   it is
reasonable to consider $k_1=k_2=:k\indexS$ and  the
interfacial temperature can be eliminated. This conceptually
reflects the modelling assumption about zero thickness of the adhesive.
Then the heat-transfer problem \eqref{heat-GC} degenerates to the transient
conditions:
\begin{align*}\\
& \frac12\big(({\bbK\indexV}_1({\thetaB}_1)\nabla{\thetaB}_1)|_{\GC}^{}\!
+({\bbK\indexV}_2({\thetaB}_2)\nabla{\thetaB}_2)|_{\GC}^{}\big){\cdot}\nu
+k\indexS\big(\JUMP{u}{\rmn},\alpha
\big)
\big({\thetaB}_1|_{\GC}{-}{\thetaB}_2|_{\GC}\big)=0,
\\ &\nonumber
\JUMP{\bbK\indexV(\thetaB)\nabla\thetaB}{\rmn}
=
\COL{-\friction(\alpha,\thetaS)\gammaCprime\big({\JUMP{u}{}}_\rmn\big)\big|\JUMP{\DT u}{\rmt}\big|-\sigma_{\rm y}(\alpha,\thetaS)\big|\DT\pi\big|}
-\big(\partial_{\DT\alpha}a_1(\JUMP{u}{},\COLOR{\alpha,}\thetaS,\DT\alpha)
\COL{-\thetaS b_0'(\alpha)}\big)\DT\alpha
\end{align*}
on $[0,T]\times\GC$, which reflects that the heat generated by delamination
is distributed  with equal proportions $\frac12$ and $\frac12$ into the two
subdomains adjacent to $\GC$. Cf.\ also \cite[Rem.\,5.3.23]{MieRou15RIS}.
Such a model has been considered e.g.\ in \cite{RosRou11TARI}.
\COL{This enhances the concept of a so-called weakly 
conducting (Kapitza) interface, cf.\ \cite{JaMBSt13TSLD}, which would be
obtained for $\friction=a_1=\sigma_{\rm y}=b_0=0$ when the normal heat flux 
would be continuous across $\GC$.}
\end{remark}





\COL{
\begin{remark}[Models with only interfacial temperature 
$\thetaS$.]\label{rem-...}\upshape
In contrast to Remark~\ref{rem-no-thetaS},
sometimes an idea of ‘an external constant temperature reservoir’ from 
\cite{BieTel01MCPF} is a relevant approximation -- then 
$\thetaB$ is considered as a fixed constant and the 
heat-transfer problem \eqref{eq6:adhes-class-form1bis}--(\ref{BC}c,d) is 
omitted, cf.\ also \cite{Roub14NRSD}.
\end{remark}
}

\REF{
\begin{remark}[Role of the plastic slip $\pi$.]\label{rem-pi}\upshape
The plastic slip $\pi$ may play various roles:\\
\Item[0.04]{1.}{To distinguish
delamination in Mode I and Mode II under gradually increasing
load as in \cite{PaMaRo??TACM,RoMaPa13QMMD}
and, after the (irreversible) delamination process is completed, it will
stop evolving,
or}
\Item[0.04]{2.}{To facilitate healing (=reversible delamination) and,
allowing for free tangential shift in the delaminated state,
to make sure that the healing will aim the shifted configuration
instead of the original one, cf.\ \cite{RoSoVo13MBMF}.}
\Item[0.04]{3.}{If $\JUMP{u}{\rmt}\sim\pi$ can be expected during the whole 
process, an approximation of friction can be realized through the activation 
threshold for evolving $\pi$.}
\end{remark}
}

\section{\COL{A monolithic form and} thermodynamics of the model
}
\label{sec-thermodyn}
Let us briefly present the thermodynamics of the boundary-value
problem \eqref{system}--\eqref{BC}. \COL{This justifies} the problem
\COL{physically. Moreover, adapting \cite[Sect.\,5.3]{MieRou15RIS}
for our situation, a ``monolithic'' 
(in the sense of simultaneous treatment of the $(d{-}1)$-dimensional 
interface and $d$-dimensional bulk) description will show a (relatively) 
simple structure of the problem, cf.\ \eqref{system-abstract} below. 

First,} the underlying overall Helmholtz {\it free energy}
$\varPsi=\varPsi(e,\REFF{\JUMP{u}{}},\pi,\alpha,\thetaS,\thetaB)
$ has a \COLOR{bulk (=volume)} and an interface parts, i.e.
\begin{subequations}\label{Psi}
\begin{align}\label{Phi-epos}
\varPsi(e,\REFF{\JUMP{u}{}},\pi,\alpha,\thetaS,\thetaB)=
\int_{\GC}\!\!\psi\indexS(\REFF{\JUMP{u}{}},
\pi,\alpha,\thetaS\COL{,\nablaS\pi,\nablaS\alpha})\,\d S
+\int_{\varOmega {\setminus} \GC}\!\!
\psi\indexV(e,\thetaB)\, \d x\,,
\end{align}
with $\psi\indexV$ and $\psi\indexS$
being the bulk and the contact-interface contributions to the specific
Helmholtz energy, respectively. One can identify
\begin{align}\nonumber
\psi\indexV (e,\thetaB)&=
\frac12\bbC(e{-}\bbE\thetaB\COLOR{{+}\bbE\theta_{_{\rm R}}}){:}
(e{-}\bbE\thetaB\COLOR{{+}\bbE\theta_{_{\rm R}}})
-\frac{(\thetaB\COLOR{{-}\theta_{_{\rm R}})}^2}2\bbB{:}\bbE-{\psi\indexV}_0(\thetaB)
\\\label{psi-bulk}
&=\frac12\bbC e{:} e-(\thetaB{-}\theta_{_{\rm R}})\bbB{:}e-{\psi\indexV}_0(\thetaB)
\quad\text{ with }\ \bbB:=\bbC\bbE.
\end{align}
Here $\frac12\bbC e{:} e$ is the {\it mechanical part of the
internal  energy in the \COLOR{bulk}}, while $-{\psi\indexV}_0(\thetaB)$ is the {\it
thermal part of the free energy}. Hereafter, we shall assume that
\begin{equation}
\label{psizero} \text{${\psi\indexV}_0: (0,+\infty)\to\R$ a strictly convex
function.}
\end{equation}
The specific \emph{contact interface} energy
$\psi\indexS(\REFF{\JUMP{u}{}},\pi,\alpha,\thetaS\COL{,\nablaS\pi,\nablaS\alpha})$ is then
\begin{equation}
\label{psi-surf}\psi\indexS\big(\REFF{\JUMP{u}{}},
\pi,\alpha,\thetaS\COL{,\nablaS\pi,\nablaS\alpha}\big)
=\begin{cases}\displaystyle{
\gamma_{\rmt}\big(\alpha,\REFF{\JUMP{u}{\rmt}}
\!{-}\pi\big)
+\gamma_{\rmn}\big(\alpha,\REFF{\JUMP{u}{\rmn}}\big)+\gammaC\big(\REFF{\JUMP{u}{\rmn}}\big)
}\hspace*{-7em}
\\[.3em]\displaystyle{\ \ \ \ \
+\frac{\kappa_{\rm H}}2|\pi|^2+\frac{\kappa_1}2|\nablaS\pi|^2+\frac{\kappa_2}2|\nablaS\alpha|^2
}\hspace*{-7em}
\\[.3em]\displaystyle{\ \ \ \ \
-\,a_0(\alpha)-b_0(\alpha)\thetaS-{\psi\indexS}_{\,0}(\thetaS)
}
&\mbox{if $\ 0\!\le\!\alpha\!\le\!1$ a.e.~on $\GC$},\\
+\infty&\text{otherwise}.\end{cases}
\end{equation}
\end{subequations}
The terms $a_0(\alpha)+b_0(\alpha)\thetaS$ represent a (temperature-dependent)
energy related to debonded adhesive, based on the phenomenon that
every damage microscopically means more microcracks/voids which always
bears some additional energy.

The typical 
\REFF{choice} 
for $\gamma_{\rmn}$ and $\gamma_{\rmt}$
describing a linearly responding (for a fixed delamination 
\REFF{variable} $\alpha$) \COLOR{elastic} adhesive
\begin{equation}\label{kappantdefs}
\gamma_{\rmn}\big(\alpha,\mbox{$\REFF{\JUMP{u}{\rmn}}$}\big)=\frac12\kappa_{\rmn}(\alpha)\big|\mbox{$\REFF{\JUMP{u}{\rmn}}$}\big|^2 \quad
\text{and} \quad \gamma_{\rmt}\big(\alpha,\mbox{$\REFF{\JUMP{u}{\rmt}}$}\big)=\frac12\kappa_{\rmt}(\alpha)\big|\mbox{$\REFF{\JUMP{u}{\rmt}}$}\big|^2,
\end{equation}
while the typical 
\REFF{choice} for $\gammaC$ is a unilateral 
normal-compliance: 
\begin{align}\label{normal-compliance}
&\gammaC\big(\REFF{\JUMP{u}{\rmn}}\big)
=\begin{cases}0&\text{if }\REFF{\JUMP{u}{\rmn}}\ge0,\\
\frac1p\kappa_{_{\rm C}}^{}\big(-\REFF{\JUMP{u}{\rmn}}\big)^p&
\text{if }\REFF{\JUMP{u}{\rmn}}<0,\end{cases}&&
\text{(unilateral normal-compliance)}
%
\intertext{which is to approximate for $\kappa_{_{\rm C}}^{}\to+\infty$ 
the Signorini unilateral contact, i.e.\ $\gammaC(\REFF{\JUMP{u}{\rmn}})=+\infty$
if $\REFF{\JUMP{u}{\rmn}}<0$ while again $\gammaC(\REFF{\JUMP{u}{\rmn}})=0$ 
otherwise, which would lead to additional variational inequality on $\GC$.
In \eqref{normal-compliance}, the exponent $p$ can be arbitrarily
large if $d=2$ or $p\le4$ if $d=3$ to be compatible with the
coercivity of the energy in the bulk (which ensures
$u\in H^1(\varOmega{\setminus}\GC;\R^d)$ so that
$\JUMP{u}{\rmn}\in H^{1/2}(\GC)\subset L^4(\GC)$
if $d=3$). 
Some situations justify also
a simpler variant of the Signorini contact, namely}
&\gammaC\big(\REFF{\JUMP{u}{\rmn}}\big)
=\begin{cases}0&\text{if }\REFF{\JUMP{u}{\rmn}}=0,\\
\infty&\text{if }\REFF{\JUMP{u}{\rmn}}\neq 0.\end{cases}&&
\text{(bilateral contact)}
\label{bilateral}\end{align}
\COL{Yet, it should be mention that} \REF{the normal compliance is more 
natural is our a setting} than \REF{the Signorini condition} which is a bit
\REF{artificial because, when the adhesive is rigid, the issue is of fracture 
and not debonding.}

The other underlying ingredient of the model is
the overall
{\it pseudo-potential of the dissipative forces} $\varXi$
which also has a bulk and
a 
\COLOR{contact-interface} contributions $\xi\indexS$ and $\xi\indexV$, namely:
\begin{align}\label{dissip-rate}
\varXi\big(\JUMP{u}{},\alpha,\thetaS,\thetaB;\DT{e},
\COLOR{\JUMP{\DT u}{}},\DT\pi,\DT{\alpha}\big)
:=
\int_{\GC}\!\xi\indexS\big(\JUMP{u}{},\alpha,\thetaS,
\COLOR{\JUMP{\DT u}{}},\DT\pi,\DT{\alpha}\big)\,\d S
+\int_{\varOmega{\setminus}\GC}\!\!\xi\indexV(\thetaB;\DT{e})\,\d x;
\end{align}
where
\begin{align*}
&\nonumber
\xi\indexS\big(\REFF{\JUMP{u}{\rmn}},\alpha,\thetaS;
\REFF{\JUMP{\DT u}{}},\DT{\alpha},\DT\pi\big)
:=
\friction(\alpha,\thetaS)
\gammaCprime(\REFF{\JUMP{u}{\rmn}})
\big|
\REFF{\JUMP{\DT u}{\rmt}}\big|+\sigma_{\rm y}(\alpha,\thetaS)\big|\DT\pi\big|
\\&\hspace{11em}
+\COL{\DT\alpha\partial_{\DT\alpha}}a_1\big(\REFF{\JUMP{\DT u}{}},\COLOR{\alpha,}\thetaS,\DT{\alpha}\big)
+\frac{d_\rmt(\alpha,\thetaS)}{\COL{2}}
\big|\REFF{\JUMP{\DT u}{\rmt}}{-}\DT\pi\big|^2
+\frac{d_\rmn(\alpha,\thetaS)}{\COL{2}}\REFF{\JUMP{\DT u}{\rmn}}^{\!\!\!\!2}\,,
\\&\xi\indexV(\thetaB;\DT{e})
:=\frac12 \bbD(\thetaB)\DT{e} {:} \DT{e}, \qquad
\end{align*}
the $a_1$-term representing 
 the specific dissipation rate
of the 
delamination process on
the contact boundary $\GC$ while the quadratic terms describes viscosity
of the adhesive. Furthermore, the overall {\it kinetic energy}
is the quadratic functional of bulk velocity \COLOR{$\DT u$, i.e.}
\begin{align}
M(\DT u)=\int_{\varOmega\setminus\GC}\frac\varrho2|\DT u|^2\d x
\end{align}
\COLOR{with $\varrho>0$ the mass density.}
Eventually, we still need the (linear) functional of external mechanical
loading $F(t)$ which, after the standard transformation to time-constant
(homogeneous) Dirichlet boundary condition by the shift $u\to u+\barudir$ with
$\barudir$ being an extension of $\udir$ from \eqref{eq6:adhes-class-form2} on
$\varOmega$, is given by
\begin{align}\label{def-of-F}
\big\langle F(t),u\rangle=\int_{\varOmega\COLOR{\setminus\GC}}\!\!\!
\big(g\COLOR{(t)}-\varrho\DDTbarudir\COLOR{(t)}\big){\cdot}u-
\big(\bbD e(\DTbarudir\COLOR{(t)})+\bbC e(\barudir\COLOR{(t)})\big){:}e(u)\,\d x
+\int_{\Gnew}\!\!f\COLOR{(t)}{\cdot}u\,\d S.
\end{align}
For simplicity and without 
{excluding} interesting applications,
we assume that $\Gdir$ and $\GC$ are far from each other so that
one can assume $\barudir|_{\SC}=0$ and thus this shift transformation
does not influence the flow rules on $\GC$.

In terms of these functionals, using also the linear operator
\begin{align}\label{def-of-E}
\boldsymbol{E}:u\mapsto\big(e(u),\JUMP{u}{}\big),
\end{align}
the mechanical part can be written
 as: 
\begin{subequations}\label{system-abstract}\begin{align}
\hspace*{-1em}M'\DDT u+E^*\partial_{(\DT{e},
\REFF{\JUMP{\DT u}{}}
)}
\varXi\big(\JUMP{u}{},\alpha,
\bftheta;\boldsymbol{E}\DT u
\big)
+\boldsymbol{E}^*\partial_{(e,\REFF{\JUMP{u}{}})}\varPsi(\boldsymbol{E}u,\pi,\alpha,
\bftheta)&\ni F(t),\ \ \ \ \
\bftheta:=(\thetaS,\thetaB),
\\
\partial_{\DT\pi}
\varXi(\JUMP{u}{},\alpha,\thetaS;
\DT\pi)
+\partial_{\pi}\varPsi\big(\JUMP{u}{},\pi,\alpha
\big)&\ni0,
\\
\partial_{\DT{\alpha}}
\varXi(\JUMP{u}{},\alpha,\thetaS;
\DT\alpha)
+\partial_{\alpha}\varPsi\big(\JUMP{u}{},\pi,\alpha,
\thetaS\big)&\ni0,
\intertext{where $\boldsymbol{E}^*$ denotes the adjoint operator 
\COLOR{to $\boldsymbol{E}$ defined in \eqref{def-of-E}} and where we \COLOR{have}
reflected a special form of $\varPsi$ from \eqref{Psi} for which 
$\partial_{\alpha}\varPsi$
is independent of $\thetaB$ and $\partial_{\pi}\varPsi$ is independent of $\thetaB$
and $\thetaS$, and also that $\partial_{(\DT{e},\REFF{\JUMP{\DT u}{}})}\varXi$ does not depend
on $\DT\pi$ and $\DT\alpha$, and similarly $\partial_{\DT\pi}\varXi$ and
$\partial_{\DT\alpha}\varXi$ do not depend on $(\DT{e},\DT v,\DT\alpha)$ and
$(\DT{e},\DT v,\DT\pi)$, respectively. This system is to be completed by
the %
{governing equations of the} heat-transfer problem \eqref{heat-GC}--\eqref{eq6:adhes-class-form1bis}--(\ref{BC}c,d)
which can be written
on such abstract level
in a ``monolithic'' form as:}
\label{abstract-heat}
\DT\bfvartheta+\mathcal{G}^*\big(\mathcal{K}(\bftheta)\mathcal{G}\bftheta\big)={\boldsymbol r}+{\boldsymbol a}
+{\boldsymbol h}_{\rm ext}\ \ \ \ \text{ with }
\ \ \ \bfvartheta=\mathcal{C}(\bftheta),&
\end{align}
\end{subequations}
where we have used the notation for
an abstract ``gradient/difference'' operator
and a heat-conductivity operator defined respectively as
\begin{subequations}\label{notation-abstract}\begin{align}\label{def-of-G}
&\mathcal{G}\bftheta:=\big(\nablaS\thetaS,\thetaB|_{\GC}\!{-}\thetaS\one,\nabla\thetaB\big)
\ \ \ \text{ and }\ \ \
\mathcal{K}(\bftheta):=\big(\bbK\indexS(\thetaS),\K\indexS(\REFF{\JUMP{u}{\rmn}},\alpha,\thetaS),
\bbK\indexV(\thetaB)\big)
\intertext{while the \COL{heat part of the internal energy, i.e.\ in fact
a} re-scaled temperature,
$\bfvartheta=(\vartheta\indexS,\vartheta\indexV)$ with
$\vartheta\indexS=C\indexS(\thetaS)$ and $\vartheta\indexV=C\indexV(\thetaB)$,
i.e.\ the heat content, is given by}
&\mathcal{C}(\bftheta):=\big(C\indexS(\thetaS),C\indexV(\thetaB)\big)
\ \ \ \text{ with }\ \ \ 
C\indexS(\thetaS):=\!\int_0^{\thetaS}\!\!\!c\indexS(t)\,\d t
\ \text{ and }\ C\indexV(\thetaB):=
\!\int_0^{\thetaB}\!\!\!c\indexV(t)\,\d t
\label{def-of-C}
\intertext{and the dissipation heat, when written very formally, is
${\boldsymbol r}=\partial_{(\DT{e},\REFF{\JUMP{\DT u}{}},\DT\alpha,\DT\pi)}
\varXi(\COLOR{\JUMP{u}{}},\alpha,\thetaS;\boldsymbol{E}\DT u,\DT\pi,\DT\alpha)\bulet(\boldsymbol{E}\DT u,\DT\pi,\DT\alpha)$,
with the ``product'' ``$\bulet$'' meant locally or, when written more
in detail component-wise, rather
}
&\nonumber
{\boldsymbol r}=\Big(
\partial_{\REFF{\JUMP{\DT{u}}{}}}\xi\indexS
\big(\JUMP{u}{},\alpha,\thetaS;\JUMP{\DT u}{\rmt}\big){\cdot}\JUMP{\DT u}{\rmt}
+\partial_{\DT\alpha}\xi\indexS\big(\JUMP{u}{},\alpha,\thetaS;
\DT{\alpha}\big)\DT\alpha
+\partial_{\DT\pi}\xi\indexS\big(\COLOR{\JUMP{u}{}},\alpha,\thetaS;
\DT\pi\big){\cdot}\DT\pi
\,,\,
\\&\hspace*{26em}
\partial_{\DT{e}}\xi\indexV\big(\thetaB;
e(\DT{u})
\big){:}e(\DT u)\Big)
\intertext{and similarly the
adiabatic contribution written component-wise is}
&{\boldsymbol a}:=\Big(
\thetaS\partial_{\alpha\thetaS}^2\psi\indexS
\big(\alpha,\thetaS\big)
\DT\alpha\,,\,
\thetaB\partial_{e\thetaB}^2\psi\indexV(e(u),\thetaB){:}e(\DT u)\Big)
=\big(-\thetaS b_0(\alpha)\DT\alpha\,,\,-\thetaB\bbB{:}e(\DT u)\big),
\intertext{and eventually the boundary heat-source from
\eqref{eq6:adhes-class-form3} is incorporated into}
&{\boldsymbol h}_{\rm ext}:=
{\big(0\,,\,h_{\rm ext}\delta_{\varGamma}^{}\big),}
\end{align}
\end{subequations}
where $\delta_{\varGamma}^{}$ denotes the surface Lebesgue measure on $\varGamma$.

\COLOR{The inclusions in (\ref{system-abstract}a,b,c) are related with 
nonsmoothness of $\varXi$ if $\friction>0$ and if $\sigma_{\rm y}>0$, and if 
$a_1(\REFF{\JUMP{u}{}},\alpha,\thetaS,\cdot)$ is nonsmooth, respectively.}
The abstract heat-transfer equation \eqref{abstract-heat} in the notation
\eqref{notation-abstract} is to be understood in the weak form as:
\begin{align}
\COLOR{\int_0^T\!\!\big\langle\hspace*{-.3em}\big\langle \mathcal{K}(\bftheta)\mathcal{G}\bftheta,\mathcal{G} {\boldsymbol v}\big\rangle\hspace*{-.3em}\big\rangle
-\big\langle{\boldsymbol r}+{\boldsymbol a}
+{\boldsymbol h}_{\rm ext},{\boldsymbol v}\big\rangle
-\big\langle\bfvartheta,\DT{\boldsymbol v}\big\rangle\,\d t=
\big\langle\bfvartheta_0,{\boldsymbol v}(0)\big\rangle}
\label{heat-eq-abstract}\end{align}
for any ${\boldsymbol v}=(v\indexS,v\indexV)\in
W^{1,\infty}(I{\times}\GC)\times W^{1,\infty}(I{\times}(\varOmega{\setminus}\GC))
$
\COLOR{with ${\boldsymbol v}(T)=0$, where the bilinear form 
$\langle\hspace*{-.25em}\langle\cdot,\cdot\rangle\hspace*{-.25em}\rangle$ is defined by}
\begin{align}\label{<<.>>}
\langle\hspace*{-.25em}\langle g,\tilde g\rangle\hspace*{-.25em}\rangle
:=\int_{\GC}g_1{\cdot}\tilde g_1+g_2{\cdot}\tilde g_2\,\d S
+\int_{\varOmega\setminus\GC}g_3{\cdot}\tilde g_3\,\d x
\end{align}
with $g=(g_1,g_2,g_3)$ and $\tilde g=(\tilde g_1,\tilde g_2,\tilde g_3)$, and
$\langle g,\tilde g\rangle:=\int_{\GC}g\indexS{\cdot}\tilde g\indexS\,\d S
+\int_{\varOmega\setminus\GC}g\indexV{\cdot}\tilde g\indexV\,\d x
+$ with $g=(g\indexS,g\indexV)$ and
$\tilde g=(\tilde g\indexS,\tilde g\indexV)$.


Introducing the \COLOR{specific} entropy ${\bfeta}=(\eta,\etaS):=-\varPsi_{\bftheta}'$, we can
write the abstract entropy equation in the form:
\begin{align}\label{entropy-eq-abstract}
\bftheta\bulet\DT\bfeta={\boldsymbol r}-\mathcal{G}^*\big(\mathcal{K}(\bftheta)\mathcal{G}\bftheta\big)
-{\boldsymbol h}_{\rm ext}.
\end{align}
This means two coupled equations for the interfacial entropy $\etaS$
and the bulk entropy $\etaB$:
\begin{subequations}\label{entropy-eq}\begin{align}\nonumber
&\thetaS\DT\eta\indexS=r\indexS+\divS(\bbK\indexS(\thetaS)\nablaS\thetaS)
+\K\indexS(\JUMP{u}{\rmn},\alpha,\thetaS)\cdot(\thetaB|_{\GC}\!{-}\thetaS\one)
\ \ \text{ with}
\\&\nonumber
\qquad\quad\
r\indexS=
\friction(\alpha,\thetaS)\gammaCprime({\JUMP{u}{}}_\rmn)|\JUMP{\DT u}{\rmt}|
+\sigma_{\rm y}(\alpha,\thetaS)|\DT\pi|
\\&\label{entropy-eq-adhes}
\qquad\qquad\
+\DT\alpha\partial_{\DT\alpha}a_1(\JUMP{u}{},\COLOR{\alpha,}\thetaS,\DT\alpha)
\COL{+d_\rmt(\alpha,\thetaS)
\big|\REFF{\JUMP{\DT u}{\rmt}}{-}\DT\pi\big|^2
+d_\rmn(\alpha,\thetaS)\REFF{\JUMP{\DT u}{\rmn}}^{\!\!\!\!2}}\ \ \text{ and}
\\\label{entropy-eq-bulk}
&\thetaB\DT\etaB=r\indexV+{\rm div}(\bbK\indexV(\thetaB)\nabla\thetaB)\ \ \text{ with}\ \
r\indexV=\bbD e(\DT{u}){:} e(\DT{u});
\end{align}\end{subequations}
here $r\indexS$ and $r\indexV$ are the heat\COL{-production rates}
in the adhesive and
\MARGINOTE{TO ${\boldsymbol r}=(r\indexS,r\indexV)$ PREDTIM NEBYLO DOBRE KDYZ $h_{\rm ext}$ BYLO ZAHRNUTO DO ${\boldsymbol r}$} in the bulk
due to \COL{the} dissipative process, 
i.e.\ ${\boldsymbol r}=(r\indexS,r\indexV)$.
Considering also the boundary conditions
\eqref{eq6:adhes-class-form3} reflected in ${\boldsymbol h}_{\rm ext}$
and \eqref{BC-heat-GC}, from \eqref{entropy-eq-bulk} we obtain
the overall entropy balance in the bulk \REFF{as} \MARGINOTE{A TEN POSLEDNI INTEGRAL V \eqref{ent-bulk} BYL PREDTIM BLBE!}
\begin{align}\nonumber
\frac{\d}{\d t}\int_{\varOmega\setminus\GC}\!\!\!\!\etaB\,\d x
&=\int_{\varOmega\setminus\GC}\!\!\!\!\frac{r\indexV
+{\rm div}(\bbK\indexV(\thetaB)\nabla\thetaB)}\thetaB\,\d x
\\\nonumber&=\int_{\varOmega\setminus\GC}\!\frac{r\indexV}\thetaB+
\frac{\bbK\indexV(\thetaB)\nabla\thetaB{\cdot}\nabla\thetaB}{{\thetaB}^2}\,\d x
\\&\qquad\ \ +\int_\varGamma\!\frac{h_{\rm ext}}\thetaB\,\d S
+\int_{\GC}\!\!\K\indexS(\JUMP{u}{\rmn},\alpha,\thetaS){\cdot}
\frac{\thetaB|_{\GC}\!{-}\thetaS\one
}{\thetaB|_{\GC}}\,\d S
\label{ent-bulk}\end{align}
where the last integral means $\int_{\GC}\sum_{i=1,2}k_i({\thetaB}_i-\thetaS)/{\thetaB}_i\,\d S$ with ${\thetaB}_i$ again meaning the traces of temperature from the
two regions adjacent to $\GC$ at a spot in question,
while from \eqref{entropy-eq-adhes} we obtain the overall entropy balance in
the adhesive
\begin{align}\nonumber
\frac{\d}{\d t}\int_{\GC}\!\!\etaS\,\d S&=\int_{\GC}\!\!\frac{r\indexS+\divS(\bbK\indexS(\thetaS)\nablaS\thetaS)
+\K\indexS(\JUMP{u}{\rmn},\alpha,\thetaS){\cdot}
(\thetaB|_{\GC}\!{-}\thetaS\one)
}{\thetaS}
\,\d S
\\&=\int_{\GC}\!\frac{r\indexS}{\thetaS}+
\frac{\bbK\indexS(\thetaS)\nablaS\thetaS{\cdot}\nablaS\thetaS}{{\thetaS}^{\!\!2}}
+\frac{\K\indexS(\JUMP{u}{\rmn},\alpha,\thetaS){\cdot}
(\thetaB|_{\GC}\!{-}\thetaS\one)
}{\thetaS}\,\d S.
\label{ent-adh}\end{align}
Altogether, summing \eqref{ent-bulk} and \eqref{ent-adh} and using the
elementary algebra $\frac{\thetaS-\thetaB}\thetaB
+\frac{\thetaB-\thetaS}{\thetaS}=\frac{(\thetaS-\thetaB)^2}{\thetaS\thetaB}$,
we obtain the {\it overall entropy balance}
\begin{align}\nonumber\\\nonumber
&\frac{\d}{\d t}\bigg(\int_{\varOmega\setminus\GC}\!\!\!\!\etaB\,\d x
+\int_{\GC}\!\!\etaS\,\d S\bigg)
=\int_{\varOmega\setminus\GC}\!\!\!\!\!\!\frac{\bbK\indexV(\thetaB)\nabla\thetaB{\cdot}\nabla\thetaB}{{\thetaB}^2}+\frac{r}\thetaB\,\d x+\int_\varGamma\!\frac{h_{\rm ext}}\thetaB\,\d S
\\[-.3em]&\qquad\qquad\qquad\qquad
+\int_{\GC}\!\!\frac{\bbK\indexS(\thetaS)\nablaS\thetaS{\cdot}\nablaS\thetaS}{{\thetaS}^{\!\!2}}
+\sum_{i=1,2}\frac{k_i\big(\JUMP{u}{\rmn},\alpha,\thetaS\big)
|{\thetaB}_i\!{-}\thetaS|^2}{\thetaS{\thetaB}_i}
+\frac{r\indexS}{\thetaS}\,\d S.
\label{entropy-overall}
\end{align}
Considering thermally isolated system (i.e.\ $h_{\rm ext}=0$), we can see that
the overall entropy $\int_{\varOmega\setminus\GC}\!\etaB\,\d x+\int_{\GC}\!\etaS\,\d S$
is nondecreasing in time, which represents the 2nd law of thermodynamics. Of
course, we relied on positive-definiteness of $\bbK\indexV$ and $\bbK\indexS$ and
non-negativity of $k\indexS$, cf.\ \eqref{ass} below. Furthermore, 
{we relied} on
non-negativity (or rather positivity) of temperatures $\thetaB$ and $\thetaS$,
which can be ensured by suitable initial/boundary conditions in 
\REFF{combination}
with the adiabatic terms $b_0'(\alpha)\thetaS$ and $\thetaB\bbC\bbE{:}e(\DT u)$
in \eqref{adhes-form-d6} and  \eqref{eq6:adhes-class-form1bis},
which may alternate sign (and hence cause both heating and cooling) but \COLOR{such possible cooling} is 
``switched off'' if the temperature $\thetaS$ approaches zero.
In other words, our system is consistent also with the 3rd law of
thermodynamics.

The mere mechanical energy balance can be obtained by testing
(\ref{system-abstract}a,b,c) respectively by $\DT u$, $\DT\pi$, and
$\DT\alpha$. We first realize the special form of 
\begin{align}
\varPsi(\boldsymbol{E}u,\pi,\alpha,\bftheta)=\mathcal{E}(\boldsymbol{E}u,\pi,\alpha)
-\varPsi_\text{\sc th}(\bftheta)+
\COL{\langle\boldsymbol b}(\boldsymbol{E}u,\pi,\alpha),\bftheta\COL{\rangle}\,,
\label{special-ansatz}\end{align}
\COL{where $\mathcal{E}=\varPsi|_{\bftheta=\bf{0}}$ is purely mechanical part
of the free energy}
with $\boldsymbol{E}$ defined in \eqref{def-of-E}, 
the linear functional
\COL{${\boldsymbol b}$}$(\boldsymbol{E}u,\pi,\alpha):\bftheta
=(\thetaS,\thetaB)\mapsto-\int_{\GC}b_0(\alpha)\thetaS\,\d S
-\int_{\varOmega{\setminus}\GC}\thetaB\bbB{:}e\,\d x$
\COL{facilitates the thermo-mechanical coupling, and where
the purely thermal part is} $\varPsi_\text{\sc th}(\bftheta):=
\int_{\GC}{\psi\indexS}_{\,0}(\thetaS)\,\d S
+\int_{\varOmega{\setminus}\GC}{\psi\indexV}_0(\thetaB)\,\d x$. The mentioned test 
then gives:
\begin{align}\nonumber
\frac{\d}{\d t}\big(M(\DT u)+\mathcal{E}(\boldsymbol{E}u,\pi,\alpha)
\big)
&+
R\big(\JUMP{u}{},\alpha,\bftheta;\boldsymbol{E}\DT u,\DT\pi,\DT\alpha)\big)
\\&=\big\langle F(t),\DT u\big\rangle
+\COL{\big\langle\partial_{e}{\boldsymbol b}(e(u))\thetaB,e(\DT u)\big\rangle
+\big\langle\partial_{\alpha}{\boldsymbol b}(\alpha)\thetaS,
\DT\alpha\big\rangle\,,}
\label{mech-energ-bal}\end{align}
where $R(\JUMP{u}{},\alpha,\bftheta;\boldsymbol{E}\DT u,\DT\pi,\DT\alpha)
=\langle{\boldsymbol r},(1,1)\rangle
=\langle\partial_{(\DT e,\REFF{\JUMP{\DT u}{}},\DT\pi,\DT\alpha)}
\varXi(\JUMP{u}{},\alpha,\bftheta;\boldsymbol{E}\DT u,\DT\pi,\DT\alpha),(\boldsymbol{E}\DT u,\DT\pi,\DT\alpha)
\rangle$
denotes the overall dissipation rate. 

Denoting by ${\boldsymbol w}=(w\indexS,w\indexV)$ the \COLOR{specific} internal energy
and by $\bfpsi=(\psi\indexS,\psi\indexV)$ the specific free energy,
we have the Gibbs' relation
${\boldsymbol w}=\bfpsi+\bftheta\bulet\bfeta
=\bfpsi|_{\bftheta=0}+\bfvartheta$.
We can write the abstract \COLOR{specific}-internal-energy balance $\DT{\boldsymbol w}
=[\partial\varPsi|_{\bftheta=\bf{0}}]\bulet(\boldsymbol{E}\DT u,\DT\pi,\DT\alpha)+{\boldsymbol r}
-\mathcal{G}^*(\mathcal{K}(\bftheta)\mathcal{G}\bftheta)$.
We denote the overall internal energy by
$$
W(\boldsymbol{E}u,\pi,\alpha,\bfvartheta)={\mathcal E}(\boldsymbol{E}u,\pi,\alpha)+H(\bfvartheta)
$$
where $H(\bfvartheta)=\langle\bfvartheta,(1,1)\rangle=\int_{\varOmega{\setminus}\GC}\!\vartheta\indexV\,\d x+\int_{\GC}\!\vartheta\indexS\,\d S$ is the total heat energy.
The {\it total-energy balance} can be obtained by testing the heat equation
by 1, i.e.\ substituting ${\boldsymbol v}=\one=(1,1)$ into
\eqref{heat-eq-abstract} and summing it with \eqref{mech-energ-bal}.
Realizing that \COLOR{$\mathcal{G}{\boldsymbol v}=0$},
it gives
\begin{align}
\frac{\d}{\d t}\big(M(\DT u)+W(\boldsymbol{E}u,\pi,\alpha,\bfvartheta)\big)
=\big\langle F(t),\DT u\big\rangle+\int_{\varGamma}h_{\rm ext}\COLOR{(t)}\,\d S.
\label{tot-energ-bal}\end{align}


\section{A semi-implicit fractional-step time discretisation}\label{sect-disc}

An efficient way for time-discretisation guaranteeing numerical stability
and convergence and, after a further spatial discretisation, allowing for
\REFF{a} constructive solution is to decouple the system by a suitable 
fractional-step method, in engineering literature usually referred to 
as a staggered scheme. \COL{Beside being computationally amenable and efficient 
conceptual algorithm,  it can simultaneously served analytically to a proof of 
a mere existence of a weak solution. For it, one can benefit from the 
monolithic formulation and use abstract results that are essentially already 
available in literature (in particular in \cite{KruRou18MMCM,MieRou15RIS})
with possibly a little modification, combining the fully-implicit
time-discretisation usually used for purely analytical purposes with some
other \cite{RouPan17QECTD}, more suitable to capture dynamical effect, as used 
also below.} \REF{In the isothermal case, also other abstract-operator-based 
techniques available for dynamic 
frictional normal-compliance contact between visco-elastic bodies exist,
cf.\ \cite{BaBaKa15DVCP,KutShi99SVMP,KutShi02DCNC} and likely they can be 
expanded if the tangential slip $\pi$ is involved, too.}

The general strategy for splitting the variables is to
\COL{respect} \REFF{the assumed} separate convexity (or at least semi-convexity)
of the free energy $\varPsi$ and additive splitting of the dissipation potential 
$\varXi$. Here, we can apply splitting of the state to $(u,\pi)$ and $\alpha$ and
$\bftheta$, assuming that 
\begin{subequations}\label{sep-convex}\begin{align}\label{sep-convex1}
&(u,\pi)\mapsto\varPsi(\boldsymbol{E}u,\pi,\alpha,\bftheta)\ \ \text{ is convex for all
$\alpha$ and $\bftheta$},
\\&\alpha\mapsto\varPsi(\boldsymbol{E}u,\pi,\alpha,\bftheta)\ \ \ \text{ is convex for all
$u$, $\pi$, and $\bftheta$},
\\&\label{sep-convex3}
c\indexS>0,\ \ \text{ and }\ \ c\indexV>0.
\end{align}\end{subequations}
\COL{In addition, we introduce a new variable $v$ in the position of 
a velocity, and in the limit indeed satisfy $v=\DT u$. This transforms 
the 2nd-order system into the 1st-order one. When combined with a 
mid-point time discretisation of the Crank-Nicolson type, it allows for 
energy-conserving scheme and for varying time step, which just for notational
simplicity is considered fixed, denoted by $\tau>0$, so that equidistant
partitions of the time interval $[0,T]$ are used. 
Let us point out that the fully implicit discretisation of the inertial term, 
i.e.\ $M'(u_\tau^k{-}2u_\tau^{k-1}{+}u_\tau^{k-2})/\tau^2$, would
serve for analytical purposes similarly good but computationally it would
lead to a scheme with spurious numerical attenuation which practically
does not allow for computation of vibrations or waves.}

This \COL{mid-point formula applied to 
the dynamical equations transformed into the form of the 1st-order system
can be understand as a particular case of the celebrated
Hilber-Hughes-Taylor formula \cite{HiHuTa77INDT}, sometimes called a 
Simo's scheme \cite[Sect.\,1.6]{SimHug98CI}. Its combination of 
fractional-step split and 1st-order-system 
transformation with mid-point formula, called sometimes the Yanenko formula 
\cite{Yane71MFS},}
leads to a three-step decoupled scheme, written in terms of
the notation from \eqref{system-abstract} as \COLOR{follows. Having,
$(u_\tau^{k-1},\alpha_\tau^{k-1},\bftheta_\tau^{k-1})$, we first solve}:
\begin{subequations}\label{system-disc}\begin{align}\nonumber
&M'
\COL{\frac{v_\tau^k{-}v_\tau^{k-1}}\tau}
+\boldsymbol{E}^*\partial_{(\DT{e},\REFF{\JUMP{\DT u}{}}
)}
\varXi\Big(\COLOR{\JUMP{u_\tau^{k-1}}{}},\alpha_\tau^{k-1},\bftheta_\tau^{k-1};
\boldsymbol{E}
\COL{v_\tau^{k-1/2}}\Big)
\\&\hspace*{9em}
+\boldsymbol{E}^*
\COL{\big[\varPsi^\circ_{(e,\REFF{\JUMP{u}{}})}}
(\COL{\,\cdot\,,\pi_\tau^{k-1/2}},\alpha_\tau^{k-1},\bftheta_\tau^{k-1})
\COL{\big](\boldsymbol{E}u_\tau^k,\boldsymbol{E}u_\tau^{k-1})}
\ni F(k\tau),\label{system-disc1}
\\
&\COL{\frac{u_\tau^k{-}u_\tau^{k-1}}\tau=\frac{v_\tau^k+v_\tau^{k-1}}2=:v_\tau^{k-1/2},}
\\
&\partial_{\DT\pi}
\varXi\Big(
{\JUMP{u_\tau^{k-1}}{},}
\alpha_\tau^{k-1},\COLOR{\thetaStaukk};
\frac{\pi_\tau^k{-}\pi_\tau^{k-1}}\tau\Big)
+\partial_{\pi}\varPsi\big(E\COL{u_\tau^{k-1/2}},\COL{\pi_\tau^{k-1/2}},
\alpha_\tau^{k-1}\big)\ni0
\intertext{\COLOR{with $\boldsymbol{E}$ defined in \eqref{def-of-E}. Second, having obtained $(u_\tau^k,\pi_\tau^k)$, we solve:}}
&\partial_{\DT{\alpha}}\varXi\Big(\COLOR{\JUMP{u_\tau^{k-1}}{}},\alpha_\tau^{k-1},\COLOR{\thetaStaukk};
\frac{\alpha_\tau^k{-}\alpha_\tau^{k-1}}\tau\Big)
+
\COL{\big[\varPsi^\circ_\alpha}
(\boldsymbol{E}u_\tau^k,\pi_\tau^k,\,\cdot\,,\COLOR{\thetaStaukk})\big]\COL{(\alpha_\tau^{k},\alpha_\tau^{k-1})}
\label{system-disc4}
\ni0.
\intertext{\COLOR{And third, having obtained also $\alpha_\tau^k$, we solve:}}
\label{heat-disc}
&\frac{\bfvartheta_\tau^k{-}\bfvartheta_\tau^{k-1}}\tau
+\mathcal{G}^*\big(\mathcal{K}(\bftheta_\tau^{k-1})\mathcal{G}\bftheta_\tau^k\big)
={\boldsymbol r}_\tau^k+{\boldsymbol a}_\tau^k+{\boldsymbol h}_{\rm ext,\tau}^k
\ \ \ \text{ with }\ \ \ \bfvartheta_\tau^k=\mathcal{C}(\bftheta_\tau^k),
\\\nonumber
&\qquad{\boldsymbol r}_\tau^k=\bigg(\partial_{
\REFF{\JUMP{\DT u}{}}}\xi\indexS\Big(\COLOR{\JUMP{u_\tau^{k-1}}{}},\alpha_\tau^{k-1},\COLOR{\thetaStaukk};
\JUMP{\COL{v^{k-1/2}}}{}
\Big){\cdot}\frac{\JUMP{
\COL{v^{k-1/2})}}{}}{\COL{1+\tau}
\epsilon^{_{(v)}}|
\JUMP{\COL{v^{k-1/2}}}{}|^2}
\\\nonumber
&\hspace{7em}
+\partial_{\DT\pi}
\xi\indexS\Big(\COLOR{\JUMP{u_\tau^{k-1}}{}},\alpha_\tau^{k-1},\COLOR{\thetaStaukk};
\frac{\pi_\tau^k{-}\pi_\tau^{k-1}}\tau\Big)
{\cdot}\frac{\pi_\tau^k{-}\pi_\tau^{k-1}}{\tau{+}\epsilon^{_{(\pi)}}|\pi_\tau^k{-}\pi_\tau^{k-1}|^2}
\\\nonumber
&\hspace{7em}
+\partial_{\DT\alpha}\xi\indexS\Big(\COLOR{\JUMP{u_\tau^{k-1}}{}},\alpha_\tau^{k-1},\COLOR{\thetaStaukk};
\frac{\alpha_\tau^k{-}\alpha_\tau^{k-1}}\tau\Big)
\frac{\alpha_\tau^k{-}\alpha_\tau^{k-1}}{\tau{+}\epsilon^{_{(\alpha)}}|\alpha_\tau^k{-}\alpha_\tau^{k-1}|^2}
\,,\,
\\\label{r-disc}
&\hspace{12em}\partial_{\DT{e}}\xi\indexV\Big(\thetaBtaukk;
e(\COL{v^{k-1/2})})
\Big){:}\frac{e
(\COL{v^{k-1/2})})}
{\COL{1+\tau}
\epsilon^{_{(e)}}|
e(\COL{v^{k-1/2})}|^2}\bigg),
\\\label{a-disc}
&\qquad{\boldsymbol a}_\tau^k:=
\bigg(
-\thetaStauk b_0(\alpha_\tau^{k-1})\frac{\alpha_\tau^k{-}\alpha_\tau^{k-1}}
{\tau+\epsilon^{(\alpha)}|\alpha_\tau^k{-}\alpha_\tau^{k-1}|^2}
\,,\,
-\thetaBtauk\bbB{:}
\frac{e
(\COL{v^{k-1/2})})}
{\COL{1+\tau}
\epsilon^{_{(e)}}|
e(\COL{v^{k-1/2})}|^2}
\bigg),
\\\label{def-of-h-tau-k}
&\qquad
{\boldsymbol h}_{\rm ext,\tau}^k:=\big(0\,,\,h_{\rm ext,\tau}^k\delta_{\varGamma}^{}\big)
\ \ \ \ \text{ with }\
h_{\rm ext,\tau}^k:=\frac1\tau\int_{(k-1)\tau}^{k\tau}
\frac{h_{\rm ext}(t,\cdot)}{1\COLOR{+\tau\epsilon^{_{(h)}}h_{\rm ext}(t,\cdot)}}\,\d t.
\end{align}
\end{subequations}
\COL{In \eqref{system-disc1} and \eqref{system-disc4}, we used the notation $(\cdot)^\circ$ for the time-difference 
defined as 
\begin{align}\label{def-of-circ}
f^\circ(z,\tilde z):=\begin{cases}(f(z)-f(\tilde z))/(z-\tilde z)&\text{if }z\ne\tilde z\\
f'(z)=f'(\tilde z)&\text{if }z=\tilde z,\end{cases}
\end{align}
assuming that $f$ is continuously differentiable. If $f(z)=\frac12z^\top Az$ is 
a quadratic form with $A$ symmetric, then 
$f^\circ(z,\tilde z)=\frac12A(z{+}\tilde z)$ gives exactly the midpoint 
difference. Cf.\ also \cite{RouPan17QECTD} for such a general discretisation 
scheme in the isothermal variant. Here, if \eqref{kappantdefs} is
adopted, it is important that the dependence of $\varPsi$ on $(\boldsymbol{E}u,\pi)$ is 
quadratic up the normal-compliance boundary term \eqref{normal-compliance}
which however depends on the scalar-valued normal displacement, so that 
the fraction in \eqref{def-of-circ} has a good sense,
and also the dependence on $\alpha$ is allowed generally nonlinear because
it is scalar-valued like the mentioned normal displacement.
For example, if the dependence on $\alpha$ would be quadratic, then 
$\varPsi^\circ_\alpha
(\boldsymbol{E}u_\tau^k,\pi_\tau^k,\,\cdot\,,\thetaStaukk)\big](\alpha_\tau^{k},\alpha_\tau^{k-1})$ in \eqref{system-disc4} would simplify to 
$\partial_\alpha\varPsi(\boldsymbol{E}u_\tau^k,\pi_\tau^k,\alpha_\tau^{k-1/2},\thetaStaukk)$.
}

Note \COL{also} that we used a regularization making ${\boldsymbol r}_\tau^k$ from
\eqref{r-disc} bounded in terms of the rates (=the time differences)
so that ${\boldsymbol r}_\tau^k$ is bounded in an $L^\infty$-space
(in contrast to ${\boldsymbol r}$
which altogether is bounded only in an $L^1$-space).
\COLOR{Similar remark concerns also $h_{\rm ext,\tau}^k$.}
 Of course, this
regularization is devised to disappear in the limit for $\tau\to0$ so that
\eqref{system-disc} has a capacity to approximate really the original
continuous system \eqref{system-abstract}, which is indeed ensured \COLOR{by the
regularizing positive coefficients $\epsilon^{_{(e)}}$, $\epsilon^{_{(v)}}$,
$\epsilon^{_{(\pi)}}$, $\epsilon^{_{(\alpha)}}$, and $\epsilon^{_{(h)}}$, which are 
of different physical dimensions. }

The decoupled system \eqref{system-disc} is to be solved recursively
for $k=1,...,T/\tau$, starting for $k=1$ with
\begin{align}\label{IC-disc}
u_\tau^0=u_0,\ \ \ \ \ \COL{v_\tau^0=}
v_0,\ \ \ \ \ \pi_\tau^0=\pi_0,
\ \ \ \ \ \alpha_\tau^0=\alpha_0,
\ \ \ \ \ \bftheta_\tau^0=
\Big(\frac{\theta_{{\scriptscriptstyle\textrm{\hspace*{-.05em}A}},0}^{}}{1{+}\tau\epsilon^{_{(h)}}
\theta_{{\scriptscriptstyle\textrm{\hspace*{-.05em}}A},0}^{}}\,,\,
\frac{{\thetaB}_{,0}}{1\COLOR{+\tau\epsilon^{_{(h)}}{\thetaB}_{,0}}}
\Big).
\end{align}

Speaking generally, desired attributes of any efficient and reliable
discretisation are the following:\\
\Item{(A1)}{existence and algorithmic amenability of solutions for the 
discrete problems,}
\Item{(A2)}{e\COL{nergetic consistency (here even energy conservation) for 
any discretisation (i.e.\ here for any the time step $\tau$),}}
\Item{(A3)}{numerical stability (=\,a-priori estimates) for any 
discretisation (i.e.\ any time step $\tau$),}
\Item{(A4)}{convergence (at least in terms of subsequences) to a suitably
defined solutions of the original continuous problem \COL{for $\tau\to0$}.}
Let us emphasize that not all these attributes are always sufficiently
justified in engineering literature.

The full model is very general and satisfaction of the above
attributes (A1)--(A4) needs relatively strong assumptions which can
be weakened in particular special cases. Considering the general
dimension $d\le3$, a rather general variant of the basic data qualification
can be 
\COL{cast} as follows:
\begin{subequations}\label{ass}
\begin{align}\label{posit}
&\bbC,\,\bbD:\R_{\mathrm{sym}}^{d\times d}\to
\R_{\mathrm{sym}}^{d\times d}\ \text{ are $4$th-order
positive definite and symmetric tensors,}
\\&\bbK\indexS{:}\R^+\to\R^{(d-1)\times(d-1)},\ \bbK\indexV{:}
\R^+\to\R^{d\times d}
\text{ continuous, 
}
\\&\label{ass-heat-capacity}
c\indexS,\,c\indexV:\R^+\to\R^+\ \text{ continuous, 
positive, with at least a linear growth
},
\\
&\COLOR{{\bbK\indexS}/{c\indexS}\,,\,{\bbK\indexV}/{c\indexV}\,\text{ bounded, uniformly positive definite},}
\\\nonumber
&\gammaC:\R\to\R^+,\COLOR{\gamma_{\rmn}:[0,1]\times\R\to\R^+,
\gamma_{\rmt}:[0,1]\times\R^{d-1}\to\R^+}
\\\label{ass-gamma-c}
&\qquad\qquad\qquad
\ \text{ continuous with less than }
\mbox{$\frac{d{+}2}{d{-}2}$}\text{-growth}, 
\text{\COLOR{and (separately) convex},}
\\&k\indexS:\R^3\to\R^+,\ \friction,d_{\rm n},d_{\rm t},\sigma_{\rm y}:\R^2\to\R^+,\
\text{ continuous, bounded,}
\\&a_0,b_0:\COLOR{[0,1]}\to\R^+\ \text{ continuous, bounded, \COLOR{and concave,}}
\\&\nonumber
a_1:\COLOR{\R^{d+3}}\to\R^+\ \text{ continuous},\ \
\COLOR{a_1(v,\alpha,\thetaS,\cdot)\text{ convex, and }}
\\&\label{ass:a1}\qquad\qquad\qquad
a_1(\JUMP{u}{},\COLOR{\alpha,}\thetaS,\cdot)\ge\epsilon|\cdot|^2
\text{ for some }\epsilon>0,
\\
&\COLOR{\friction,d_{\rmt},d_{\rmn}:[0,1]\times\R^+\to\R^+\ \text{ continuous, bounded},}
\\&\varrho\ge0,\ \kappa_{\rm H}>0,\ \kappa_1>0,\ \kappa_2>0,
\\&\barudir\!\in\!H^1(I;H^1(\varOmega{\setminus}\GC;\R^d)),
\ \varrho\barudir\!\in\!W^{2,1}(I;L^2(\varOmega;\R^d)),
\ g\!\in\!L^2(Q;\R^d),\ f\!\in\!L^2(\Snew;\R^d),
\\&u_0\!\in\!H^1(\varOmega{\setminus}\GC;\R^d),\ \ 
v_0\!\in\!L^2(\varOmega;\R^d),\ \ \pi_0\!\in\!H^1(\GC;\R^{d-1}),\ \
\alpha_0\!\in\!H^1(\GC),\ \ 0\le\alpha_0\le1,
\\&h_{\rm ext}\!\in\!L^1(\Sigma),\ \ \ \ h_{\rm ext}\ge0,\ \ \ \ 
c\indexS(\theta_{{\scriptscriptstyle\textrm{\hspace*{-.05em}A}},0}^{})\!\in\!L^1(\GC),\ \
\theta_{{\scriptscriptstyle\textrm{\hspace*{-.05em}A}},0}^{}\ge0,\ \ \ \ 
c\indexV({\thetaB}_{,0})\!\in\!L^1(\varOmega),\ \ \ \ {\thetaB}_{,0}\ge0,
\end{align}\end{subequations}
where $\barudir$ is the extension of $\udir$ used in \eqref{def-of-F}
and $\R^+$ is the set of non-negative reals. Here the symmetry \eqref{posit}
means that $\bbC_{ijkl}= \bbC_{jikl}= \bbC_{klij}$, and the same for $\bbD$.
We have used the standard notation for the function
spaces on $\varOmega$ and their norms. Namely,
$L^p(\varOmega)$
denotes the Lebesgue space of measurable functions whose $p$-power in
integrable while $W^{k,p}(\varOmega)$ denotes the Sobolev space of functions
which are together with all their $k^{\rm th}$-order derivatives in $L^p(\varOmega)$.
Moreover, we abbreviate $W^{k,2}(\varOmega)=H^k(\varOmega)$, as usual.
Below in (\ref{est}g,h), $(\cdot)^*$ refers to the dual space.
For vector or matrix valued cases, we will write e.g.\ $L^p(\varOmega;\R^n)$
or $L^p(\varOmega;\R^{n\times m})$ etc.
Moreover, we use the Bochner spaces of Banach-space-valued functions on
the time interval $I=[0,T]$, denoted by $L^p(I;X)$.

In fact, \eqref{est-DTu} requires controlled growth of the boundary
conditions, in particular \eqref{ass-gamma-c}, and excludes the 
\COLOR{mentioned} Signorini contact 
in the dynamical case (where $\varrho>0$ and thus
$M\ne0$) while the bilateral contact
\eqref{bilateral} works if the underlying space $H^1(\varOmega;\R^d)$
is replaced by the corresponding closed linear subspace respecting
\eqref{bilateral}.

\COL{The following assertion concerns the attributes (A1) and (A2). 
Considering a fixed time step $\tau>0$,
we define the piecewise-constant and the piecewise 
affine interpolants respectively by 
\begin{subequations}\label{def-of-interpolants}
\begin{align}\label{def-of-interpolants-}
&&&
\overline{u}_\tau(t)= u_\tau^k,\qquad\ \
\underline u_\tau(t)= u_\tau^{k-1},\qquad\ \
\underline{\overline u}_\tau(t)=\frac12u_\tau^k+\frac12u_\tau^{k-1},
&&\text{and}
&&
\\&&&\label{def-of-interpolants+}
u_\tau(t)=\frac{t-(k{-}1)\tau}\tau u_\tau^k
+\frac{k\tau-t}\tau u_\tau^{k-1}
&&\hspace*{-8em}\text{for }(k{-}1)\tau<t\le k\tau.
\end{align}\end{subequations}
Similar meaning has also $\alpha_\tau$, $\overline\alpha_\tau$, 
$\underline{\overline v}_\tau$, $\overline f_\tau$, etc.

\begin{proposition}[Existence and energetics of approximate solutions.]
Let \eqref{ass} holds. Then, for any $k=1,...,T/\tau$, there exists the 
solution $(u_\tau^k,v_\tau^k,\pi_\tau^k,\alpha_\tau^k,\bftheta_\tau^k)
\equiv(u_\tau^k,v_\tau^k,\pi_\tau^k,\alpha_\tau^k,\thetaStauk,\thetaBtauk)
\in H^1(\varOmega{\setminus}\GC;\R^d)^2\times H^1(\GC;\R^{d-1})\times H^1(\GC)\times 
H^1(\GC)\times H^1(\varOmega{\setminus}\GC)$ to the scheme 
\eqref{system-disc}--\eqref{IC-disc}.
Moreover, the discrete mechanical energy balance, i.e.\ an analog
of \eqref{mech-energ-bal} integrated over time, holds
\begin{align}\nonumber
M(v_\tau(t))+\mathcal{E}(\boldsymbol{E}u_\tau(t),\pi_\tau(t),\alpha_\tau(t))
+
\int_0^t\!\!
R\big(\JUMP{\underline u_\tau}{},\underline\alpha_\tau,\underline\bftheta_\tau;
\boldsymbol{E}\underline{\overline v}_\tau,\DT\pi_\tau,\DT\alpha_\tau)\big)\,\d t\ \ \ \ \ \ \ \ \ \ 
\\=\int_0^t\!\!\Big(\big\langle\overline F_\tau,\underline{\overline v}_\tau\big\rangle
+\COL{\big\langle\partial_{e}{\boldsymbol b}(e(u_\tau))\underthetaBtau,e(\underline{\overline v}_\tau)\big\rangle
+\big\langle\partial_{\alpha}{\boldsymbol b}(\underline\alpha_\tau)\underthetaStau,
\DT\alpha_\tau\big\rangle\Big)\,\d t\,,}
\label{engr-mech-disc}
\end{align}
and, at least if there is no adiabatic coupling ${\boldsymbol b}=0$,
the analog of the total energy balance\eqref{tot-energ-bal} integrated over 
time holds at least
as an inequality
\begin{align}\nonumber
M(v_\tau(t))
+W(\boldsymbol{E}u_\tau(t),\pi_\tau(t),\alpha_\tau(t),\bfvartheta_\tau(t))
\le\int_0^t\!\!\Big(\big\langle\overline F_\tau,\underline{\overline v}_\tau\big\rangle
+\int_{\varGamma}h_{\rm ext}\,\d S\Big)\,\d t\ \ \ \ \ \ \ 
\\[-.0em]+M(v_0)+W(\boldsymbol{E}u_0,\pi_0,\alpha_0,\bfvartheta_0)
\label{engr-tot-disc}
\end{align}
for any $t=k\tau$ with $k=1,...,T/\tau$. Moreover,
the non-negativity of $\bftheta_\tau^k$ on $\GC\times\varOmega$ holds, too.
\end{proposition}
}

\noindent{\it Sketch of the proof.}
As to the attribute (A1), 
note that (\ref{ass}a,b,e,g) implies \eqref{sep-convex}. 
%
Further, it is important that the scheme \eqref{system-disc} 
has a variational character 
for all three sub-problems for $(u_\tau^k,\COL{v_\tau^k},\pi_\tau^k)$, for
$\alpha_\tau^k$, and for $\bftheta_\tau^k$. These potentials are coercive and
weakly lower semicontinuous, which ensures by the standard direct method
that a solution of this system always does exist. \COL{Noteworthy, 
even though the system for the first step (\ref{system-disc}a-c) is 
nonsymmetric, it indeed has a potential after elimination of 
$v_\tau^k=2(u_\tau^k{-}u_\tau^{k-1})/\tau-v_\tau^{k-1}$, namely
\begin{align}\nonumber
(u,\pi)\mapsto 2\tau M\Big(\frac{u{-}\tau v_\tau^{k-1}{-}u_\tau^{k-1}}{\tau^2}\Big)
+\varXi\Big(\JUMP{u_\tau^{k-1}}{},\alpha_\tau^{k-1},\bftheta_\tau^{k-1};
\boldsymbol{E}\frac{u{-}u_\tau^{k-1}}\tau,\frac{\pi{-}\pi_\tau^{k-1}}\tau\Big)
\\
+\frac2\tau\varPsi\Big(\boldsymbol{E}\frac{u{+}u_\tau^{k-1}}2,
\frac{\pi{+}\pi_\tau^{k-1}}2,\alpha_\tau^{k-1},\bftheta_\tau^{k-1}\Big)
-\Big\langle F(k\tau),\frac u\tau\Big\rangle
\label{potential-for-upi}\end{align}
if $\varPsi(\cdot,\cdot,\cdot,\alpha,\bftheta)$ would be quadratic (i.e.\ 
$\gammaC=0$) while for the nonquadratic case the primitive of the difference 
quotients \eqref{def-of-circ} are to be involved, cf.\ 
\cite[Formula (4.4)]{RouPan17QECTD}.
Also,} note that, due to the regularization of the heat-sources in 
(\ref{system-disc}e-g),
we can use the conventional $L^2$-theory for the elliptic boundary-value
problems.

\COL{
Testing the particular equations (\ref{system-disc}a,b,c) respectively by 
$v_\tau^{k-1/2}$, $u_\tau^{k-1/2}$, 
and $(\pi_\tau^{k}{-}\pi_\tau^{k-1})/\tau$, we obtain

\begin{align}\nonumber
&M(v_\tau^k)+
\mathcal{E}(\boldsymbol{E}u_\tau^k,\pi_\tau^k,\alpha_\tau^{k-1})
+\tau R\Big(\JUMP{u_\tau^{k-1}}{},\alpha_\tau^{k-1},\bftheta_\tau^{k-1};
\boldsymbol{E}v_\tau^{k-1},\frac{\pi_\tau^k{-}\pi_\tau^{k-1}}\tau,
0\Big)
\\
&=M(v_\tau^{k-1})+\mathcal{E}(\boldsymbol{E}u_\tau^{k-1},\pi_\tau^{k-1},\alpha_\tau^{k-1})
+F_\tau^k\cdot v_\tau^{k-1/2}+
\tau\big\langle{\partial_e\boldsymbol b}(e(u_\tau^{k-1}))\thetaBtauk\,\,
v_\tau^{k-1/2}\big\rangle.
\label{system-disc1-3-tested}\end{align}
For this we used, among others, the binomial formulas
\begin{align}\nonumber
&\varrho\frac{v_\tau^k{-}v_\tau^{k-1}}\tau{\cdot}v_\tau^{k-1/2}
=\varrho\frac{v_\tau^k{-}v_\tau^{k-1}}\tau{\cdot}
\frac{v_\tau^k{+}v_\tau^{k-1}}2=\frac\varrho2|v_\tau^k|^2-\frac\varrho2|v_\tau^{k-1}|^2
\ \text{ and }\ 
\\&\nonumber
\bbC e(u_\tau^{k-1/2}){\colon}e(v_\tau^{k-1/2})
=\bbC e\Big(\frac{u_\tau^k{+}u_\tau^{k-1}}2\Big){\colon}e\Big(\frac{u_\tau^k{-}u_\tau^{k-1}}\tau\Big)=\frac12\bbC e(u_\tau^k){\colon}e(u_\tau^k),\ \text{ or, by
\eqref{def-of-circ}, also}\
\\&\nonumber
{\gammaC}^{\!\!\circ}(\JUMP{u_\tau^k}{\rmn},\JUMP{u_\tau^{k-1}}{\rmn})\JUMP{v^{k-1/2}}{\rmn}=
\frac{\gammaC(\JUMP{u_\tau^k}{\rmn})-\gammaC(\JUMP{u_\tau^{k-1}}{\rmn})}{\JUMP{u_\tau^k}{\rmn}
-\JUMP{u_\tau^{k-1}}{\rmn}}\JUMP{v^{k-1/2}}{\rmn}
\\&\nonumber\qquad\qquad\qquad=
\frac{\gammaC(\JUMP{u_\tau^k}{\rmn})-\gammaC(\JUMP{u_\tau^{k-1}}{\rmn})}{\JUMP{u_\tau^k}{\rmn}-\JUMP{u_\tau^{k-1}}{\rmn}}\JUMP{\frac{u_\tau^k-u_\tau^{k-1}}\tau}{\rmn}
=\frac{\gammaC(\JUMP{u_\tau^k}{\rmn})-\gammaC(\JUMP{u_\tau^{k-1}}{\rmn})}\tau
\end{align}
a.e.\ on $\varOmega$ and on $\GC$, respectively. Then, testing 
\eqref{system-disc4} by $\alpha_\tau^k{-}\alpha_\tau^{k-1}$, gives
\begin{align}\nonumber
&\mathcal{E}(\boldsymbol{E}u_\tau^k,\pi_\tau^k,\alpha_\tau^{k})
+\tau R\Big(\JUMP{u_\tau^{k-1}}{},\alpha_\tau^{k-1},\bftheta_\tau^{k-1};0,0,
\frac{\alpha_\tau^k{-}\alpha_\tau^{k-1}}\tau\Big)
\\&\qquad
=
\mathcal{E}(\boldsymbol{E}u_\tau^k,\pi_\tau^k,\alpha_\tau^{k-1})
+
\tau\big\langle\partial_\alpha{\boldsymbol b}(e(u_\tau^{k-1}))\thetaBtauk\,\,
v_\tau^{k-1/2}\big\rangle\,.
\label{system-disc4-tested}
\end{align}
Summing \eqref{system-disc1-3-tested} and \eqref{system-disc4-tested} up,
we can enjoy cancellation of the terms 
$\pm\varPsi(\boldsymbol{E}u_\tau^k,\pi_\tau^k,\alpha_\tau^{k-1},\thetaStaukk)$
and, when still sum for $k=1,2,...l\le T/\tau$, we obtain
\begin{align}\nonumber
&M(v_\tau^l)+\mathcal{E}(\boldsymbol{E}u_\tau^l,\pi_\tau^l,\alpha_\tau^{l})
+\tau\sum_{k=1}^lR\Big(\JUMP{u_\tau^{k-1}}{},\alpha_\tau^{k-1},\bftheta_\tau^{k-1};
\boldsymbol{E}v_\tau^{k-1/2},\frac{\pi_\tau^k{-}\pi_\tau^{k-1}}\tau,
\frac{\alpha_\tau^k{-}\alpha_\tau^{k-1}}\tau\Big)
\\&\nonumber=M(v_0)+\mathcal{E}(\boldsymbol{E}u_0,\pi_0,\alpha_0)
\\&\quad
+\tau\sum_{k=1}^l\big\langle F_\tau^k,v_\tau^{k-1/2}\big\rangle
+\big\langle\partial_e{\boldsymbol b}(e(u_\tau^{k-1}))\thetaBtauk\,\,
v_\tau^{k-1/2}\big\rangle
+\Big\langle\partial_{\alpha}{\boldsymbol b}(\alpha_\tau^{k-1})\thetaStauk,
\frac{\alpha_\tau^k{-}\alpha_\tau^{k-1}}\tau\Big\rangle\,.
\label{system-disc1-4-tested}
\end{align}
Thus \eqref{engr-mech-disc} has been obtained.

Next, we test the heat-transfer problem (\ref{system-disc}e-h) by $\one=(1,1)$.
When added to \eqref{system-disc1-4-tested}, this ``nearly'' cancels the 
dissipation $R$-term up to the regularizing $\epsilon$-coefficients; i.e.\
the cancellation would really occur if these coefficients would vanish, while
otherwise this regularized dissipation heat is always dominated by 
the $R$-term in \eqref{system-disc1-4-tested} and thus it gives the 
inequality in \eqref{engr-tot-disc}.

Eventually,} we have to show the non-negativity of $\bftheta_\tau^k=
(\vartheta_{{\scriptscriptstyle\textrm{\hspace*{-.05em}A}},\tau}^k,{\vartheta\indexV}_{,\tau}^k)$
a.e.\ on $(\GC\times\varOmega{\setminus}\GC)$. This can be performed
standardly by the test of the negative part of 
$\bftheta_\tau^k$, still belonging to $H^1(\GC)\times H^1(\varOmega{\setminus}\GC)$,
hence being a legitimate test function.
$\hfill\Box$

\medskip

Moreover, if \eqref{sep-convex} holds and in addition also
\begin{align}\label{adiab-convex}
&\partial_{\alpha\thetaS}^2\psi\indexS(v,\pi,\alpha,\thetaS)\ \text{ and }\ 
  \partial_{e\thetaB}^2\psi\indexV(e,\thetaB)\ \text{ are constant
    in $(\thetaS,\thetaB)$},
\end{align}
then all the underlying potentials are convex, which eliminates the usual 
algorithmic difficulties with finding critical points of nonconvex functionals 
and thus makes the numerical solution relatively simple. Note that 
\eqref{sep-convex3} ensures that $C\indexS$ and $C\indexV$, determining 
$\mathcal{C}(\cdot)$ in \eqref{heat-disc}
and defined by \eqref{def-of-C}, are monotone so that their primitive
functions occurring in the potential underlying for (\ref{system-disc}d-f) are
convex, and similarly the contribution \REFF{coming} from \eqref{a-disc}
is convex (and even quadratic) if \eqref{adiab-convex} holds.
Often, these minimization problems have (possibly after a Mosco-type
transformation \cite[Sect.\,3.6.3]{MieRou15RIS}) a structure of
linear-quadratic programming or second-order cone programming problems for
which efficient computational algorithms are widely available.
If healing is not allowed (thus always $\alpha_\tau^k{-}\alpha_\tau^{k-1}\le0$),
one can weaken \eqref{adiab-convex} by allowing
$\partial_{\alpha\thetaS}^2\psi\indexS$ depending on $\thetaS$ in a way so that
$\thetaS\mapsto\thetaS\partial_{\alpha\thetaS}^2\psi\indexS\big(\JUMP{u}{},\pi,\alpha,\thetaS\big)$ is nondecreasing, which still makes the contribution
to the underlying potential coming from  \eqref{a-disc} convex (although not
quadratic). In more general cases when \eqref{adiab-convex} is not
satisfied, the potential of (\ref{system-disc}d-f) may be nonconvex and
finding its critical point (not necessarily a global minimum) by
some iterative routine {is necessary}.

Another important attribute of the scheme \eqref{system-disc} is
$\bftheta_\tau^k\ge0$. This is facilitated by \COL{the adiabatic heat sources} 
${\boldsymbol a}_\tau^k$ involving $\bftheta_\tau^k$
instead of $\bftheta_\tau^{k-1}$ occurring in (\ref{system-disc}a,b) which,
on the other hand, needs a bit stronger qualification of the data
as far as growth conditions than physically necessary,
namely the heat capacities are required to depend on temperature with
an at least linear growth, see \eqref{ass-heat-capacity}.

As to the attribute (A3), one can imitate the procedure which led to
\eqref{mech-energ-bal} and \eqref{tot-energ-bal} when replacing
time derivatives by time differences. It allows for optimal data
qualification when rather technical 
interpolation estimates
for the adiabatic term is used, cf.\ e.g.\ \cite{RosRou11TARI,RosRou13ACDM}.
Alternatively, 
we use here a simpler \COL{based on}  \eqref{disc-balance} below,
although it needs a data qualification \COL{(\ref{ass}c,d) which is a bit 
stronger in comparison what would be needed when used the mentioned 
interpolation}.

The general physically motivated strategy for (A3) is to get the mechanical
and the heat energies $M+{\cal E}_0+H$ bounded uniformly in time and the
overall energy dissipated during the process $\int_0^T\varXi\,\d t$ bounded,
both uniformly in $\tau>0$, too. \COLOR{Based on estimation of 
$M+{\cal E}_0+\frac12H$, we have:}

\COLOR{
\begin{proposition}[Numerical stability]\label{prop-1}
Let 
\eqref{ass} hold. 
With $C$ and $C_p$ independent of $\tau>0$, 
the following a-priori estimates hold:
\begin{subequations}\label{est}\begin{align}
&\big\|u_\tau\big\|_{H^1(I;H^1(\varOmega\COLOR{\setminus\GC};\R^d))}\le C\ \ \ \text{ and }\ \ \
\big\|M\DT u_\tau\big\|_{L^{\infty}(I;L^2(\varOmega;\R^d))}\le C,
\\&\big\|\pi_\tau\big\|_{L^{\infty}(I;H^1(\GC;\R^{d-1}))\,\cap\,H^1(I;L^1(\GC))}\le C,
\\&\big\|\alpha_\tau\big\|_{L^{\infty}(\SC)\,\cap\,L^{\infty}(I;H^1(\GC))
\,\cap\,H^1(I;L^2(\GC))}\le C,
\\&\label{est-L1-theta}
\big\|\vartheta_{{\scriptscriptstyle\textrm{\hspace*{-.05em}A}},\tau}^{}
\big\|_{L^\infty(I;L^1(\GC))}\le C\ \ \ \text{ and }\ \ \
\big\|{\vartheta\indexV}_{,\tau}\big\|_{L^\infty(I;L^1(\varOmega))}\le C,
\\&\big\|\nablaS\vartheta_{{\scriptscriptstyle\textrm{\hspace*{-.05em}A}},\tau}^{}\big\|_{L^p(\SC;\R^{d-1})}\le C_p\ \ \text{ with }\ 1\le p<1+1/d,
\\
&\big\|\nabla{\vartheta\indexV}_{,\tau}\big\|_{L^p(Q;\R^d))}\le C_p\ \ \ \ \ \ \
\text{ with }\ 1\le p<(d{+}2)/(d{+}1),
\\
&\label{est-DTu}
\big\|M\DT u_\tau\big\|_{H^1(I;H^1(\varOmega\setminus\GC;\R^d)^*)}\le C,
\\&\label{est-DTtheta}
\big\|\DT\vartheta_{{\scriptscriptstyle\textrm{\hspace*{-.05em}A}},\tau}^{}
\big\|_{L^1(I;H^2(\GC)^*)}\le C\ \ \ \text{ and }\ \ \
\big\|\DT{\vartheta\indexV}_{,\tau}\big\|_{L^1(I;H^3(\varOmega\setminus\GC)^*)}\le C,
\end{align}\end{subequations}
where $u_\tau$ denotes the continuous piece-wise affine interpolant
from the values $(u_\tau^k)_{k=0}^{T/\tau}$ and similarly also $\pi_\tau$,
$\alpha_\tau$, etc. 
\end{proposition}
}

\noindent{\it Sketch of the proof.}
To show numerical stability of the scheme
is to test the particular equations/inclusions in \eqref{system-disc}
successively by \COL{$v_\tau^{k-1/2}$, $u_\tau^{k-1/2}$,
$(\pi_\tau^k{-}\pi_\tau^{k-1})/\tau$,
$(\alpha_\tau^k{-}\alpha_\tau^{k-1})/\tau$,} and $1/2$ (meaning a constant 1/2 on
both $\varOmega{\setminus}\GC$ and $\GC$) instead of $\one$ used for
\eqref{tot-energ-bal}; cf.\ \cite[Sect.5.3]{MieRou15RIS}. Under the assumption
(\ref{sep-convex}a,b), one can see a telescopic cancellation effect and the
resulted estimate: 
\begin{align}\nonumber
&
M(\COL{v_\tau^{k}})
+\mathcal{E}\big(\boldsymbol{E}u_\tau^k,\pi_\tau^k,\alpha_\tau^k\big)
+\frac12 H(\bfvartheta_\tau^k)
\\[-.7em]&\nonumber
\qquad\qquad
+\frac12\sum_{l=1}^k\varXi\Big(\COLOR{\JUMP{u_\tau^{l-1}}{}},\alpha_\tau^{l-1},\COLOR{\thetaStaull};\boldsymbol{E}
\COL{v_\tau^{k-1/2}},
\frac{\pi_\tau^l-\pi_\tau^{l-1}}\tau,
\frac{\alpha_\tau^l-\alpha_\tau^{l-1}}\tau\Big)
\\[-.7em]&\nonumber
\le
 M(v_0)+\mathcal{E}_0\big(\boldsymbol{E}u_0,\pi_0,\alpha_0)
+\frac12 H(\bfvartheta_0)
+\sum_{l=1}^k\big\langle F_\tau^l,
\COL{v_\tau^{k-1/2}}\big\rangle
+\int_{\varGamma}h_{\rm ext,\tau}^l\,\d S
\\[-.7em]&\qquad\qquad
+\sum_{l=1}^k\bigg\langle
\mathcal{E}_1(\boldsymbol{E}u_\tau^l,\pi_\tau^l,\alpha_\tau^{l-1})\Big(\bftheta_\tau^{l-1}\!
-\frac12\bftheta_\tau^l\Big),
\Big(\boldsymbol{E}
\COL{v_\tau^{k-1/2}},
\frac{\pi_\tau^l{-}\pi_\tau^{l-1}}\tau,
\frac{\alpha_\tau^l{-}\alpha_\tau^{l-1}}\tau\Big)\bigg\rangle,
\label{disc-balance}
\end{align}
where $\bfvartheta_0=\mathcal{C}(\bftheta_0)$, $F_\tau^k=F(k\tau)$, and
$h_{\rm ext,\tau}^k$ is from \eqref{def-of-h-tau-k}.
Here we have used the special ansatz leading to \eqref{special-ansatz}.
By using the discrete Gronwall
inequality for \eqref{disc-balance} and by the coercivity of the mentioned
energies, it then yields the a-priori estimates (\ref{est}a-d).

Moreover, special nonlinear tests
\COL{by $\one-\one/(\one{+}\bfvartheta)^\epsilon$ with $\epsilon>0$}
allow for important estimates of the 
temperature gradients, namely  (\ref{est}e-f), cf.\ \cite{BocGal89NEPE} 
or also e.g.\ \cite{RosRou11TARI,MieRou15RIS}. 
\COLOR{There are several variants \cite[Sect.\,8.3]{KruRou18MMCM}
how to 
\REFF{perform} this test. The simplest is to substitute
for $\bftheta=\mathcal{C}^{-1}(\bfvartheta)$ and to eliminate $\bftheta$
by this way. Then we can rely on 
\eqref{est-L1-theta} to derive (\ref{est}e-f).}


Eventually, by comparison from the equations, we obtain still the 
``dual'' estimates 
\eqref{est-DTtheta}.
$\hfill\Box$

\medskip

\REFF{Furthermore,} the \COLOR{desired} attribute (A4), 
i.e.\ convergence, is based on
several carefully assembled arguments to be applied in a proper order.
\COLOR{The definition of a weak solution is to take into account 
that (\ref{system-abstract}a,b,c) are inclusions involving convex 
subdifferentals if $\friction>0$ and if $\sigma_{\rm y}>0$ and if 
$a_1(v,\alpha,\thetaS,\cdot)$ is nonsmooth,
which results to variational inequalities rather than equalities.
More specifically, (\ref{system-abstract}a) says that 
$F(t)-M'\DDT u-E^*\partial_{(\DT{e},
\REFF{\JUMP{\DT u}{}}
)}
\varXi\big(\JUMP{u}{},\alpha,
\bftheta;\boldsymbol{E}\DT u
\big)
-E^*\partial_{(e,\REFF{\JUMP{u}{}})}\varPsi(\boldsymbol{E}u,\pi,\alpha,\bftheta)$ belongs to the subdifferential 
of the convex functional $\DT u\mapsto \varXi(\JUMP{u}{},\alpha,\bftheta;\boldsymbol{E}\DT u)$,
which straightforwardly \COL{(after integrating the term
$\langle M'\DDT u,\DT u\rangle$ by part in time)}
  leads to the weak formulation of 
(\ref{system-abstract}a) as
\begin{subequations}\label{weak-sln}\begin{align}\nonumber
&\int_0^T\!\!
\varXi\big(\JUMP{u}{},\alpha,\bftheta;E\wt u\big)
+\Big\langle\partial_{(e,\REFF{\JUMP{u}{}})}\varPsi(\boldsymbol{E}u,\pi,\alpha,\bftheta),
E(\wt u{-}\DT u)\Big\rangle
-\big\langle F,\wt u{-}\DT u\big\rangle-\big\langle M'\DDT u,\DT{\wt u}\big\rangle\,\d t
\\[-.6em]&\hspace{9em}
\ge
M(\DT u(T))+\int_0^T\!\!\varXi\big(\JUMP{u},\alpha,\bftheta;\boldsymbol{E}\DT u\big)\,\d t
+\big\langle M'v_0,\wt u(0)\big\rangle-M(v_0)
\intertext{for all $\wt u\in H^1(I{\times}(\varOmega{\setminus}\GC);\R^d)$ 
with $\wt u(T)=0$, while (\ref{system-abstract}b,c) leads to}
&\int_0^T\!\!
\partial_{\pi}\varPsi\big(\boldsymbol{E}u,\pi,\alpha,
\thetaS\big)(\wt\pi{-}\DT\pi)+
\varXi\big(\JUMP{u}{},\alpha,\thetaS;0,0,\wt\pi,0\big)\,\d t\ge
\int_0^T\!\!\varXi\big(\JUMP{u}{},\alpha,\thetaS;
0,0,\DT\pi,0\big)\,\d t,\ 
\\&
\int_0^T\!\!
\partial_{\alpha}\varPsi\big(\boldsymbol{E}u,\pi,\alpha,
\thetaS\big)
(\wt\alpha{-}\DT\alpha)+
\varXi\big(\JUMP{u},\alpha,\thetaS;
0,0,0,\wt\alpha\big)\,\d t\ge
\int_0^T\!\!\varXi\big(\JUMP{u}{},\alpha,\thetaS;
0,0,0,\DT\alpha\big)\,\d t\
\end{align}\end{subequations}
for all $\wt\pi\in L^2(I;H^1(\GC;\R^d))$ and $\wt\alpha\in L^2(I;H^1(\GC))$,
while the heat equation \eqref{abstract-heat}--\eqref{notation-abstract} 
with the initial condition from \eqref{IC} has been weakly formulated already 
as \eqref{heat-eq-abstract}, and eventually the resting initial conditions
are satisfied, i.e.\ $u(0)=u_0$, $\pi(0)=\pi_0$, and $\alpha(0)=\alpha_0$.
}

\begin{proposition}[Convergence]
Let again \eqref{ass} hold and let $(u_\tau,\pi_\tau,\alpha_\tau,
{\thetaS}_\tau,{\thetaB}_\tau)$ be some approximate 
solution from Proposition~\ref{prop-1}. Then there is a subsequence with 
$\tau\to0$ converging weakly* in the topologies from (\ref{est}a-c,e-g) to 
some limit $(u,\pi,\alpha,\thetaS,\thetaB)$.
Moreover, each 
\REFF{of these limits} is a weak solution in the sense 
\eqref{weak-sln} with \eqref{heat-eq-abstract} to the original problem
\eqref{system-abstract} with the initial conditions \eqref{IC}. 
\end{proposition}

\noindent{\it Sketch of the proof.}
%
The general strategy is, after 
\COLOR{formulating the discrete scheme \eqref{system-disc} in the form
like \eqref{weak-sln} and \eqref{heat-eq-abstract} by using the by-part
summation in time and after}
selecting weakly* convergent subsequences
by Banach's selection principle, to make the limit passage in the mechanical
part \COLOR{by using also the compactness arguments through the 
Aubin-Lions theorem}, and further to show the mechanical-energy conservation 
and to make limit passage in the total dissipated energy as an equality, which 
further shows the strong convergence of the dissipative heat. \DELETE{and} This 
eventually allows for a limit passage in the heat-transfer equation.
We refer e.g.\ to a rather abstract scheme \cite[Sect.5.3.2]{MieRou15RIS}
which covers our general model under the assumptions \eqref{ass}.
$\hfill\Box$

\medskip

\begin{remark}[Numerical implementation.]\upshape
Coming back to (A1), the numerical implementation of \eqref{system-disc}
needs still a spatial discretization. To this goal, a finite-element method
is well applicable. One can consider P1 or Q1 finite elements in bulk
for discretisation of $u$, $\thetaB$, and $\vartheta\indexV$, and their 
\COLOR{interface} variants for discretisation of $\pi$, $\alpha$, 
$\thetaS$, and $\vartheta\indexS$\COL{; cf.\ e.g.\ \cite{EllRan13FEAC} for
  coupled bulk-surface FEM}. 
In any case, in coupled anisothermal
situations, the adiabatic-like effects caused by \eqref{a-disc} makes
the non-negativity of temperature in spatially discretised problem
questionable even on acute triangulations and only a successive limit
passage seems to work, i.e.\ first the space discretisation and afterwards
the time discretisation; cf.\ also the discussion in
\cite{BarRou11TVER,BarRou13NATC}.
\end{remark}

\begin{remark}[The role of $a_1$ and $a_0$.]\upshape
\COL{The dissipation and stored interfacial energies
$a_1$ and $a_0$ should be distinguished in particular for
anisothermal or reversible-delamination situations. Typically,
$a_1$ involves rate-dependent effects reflecting the idea that fast
processes leads to some extra microscopical atomic vibrations, i.e.\
some heat production, while $a_0$ reflects the idea that degradation of
the adhesive creates microscopical new interface or microcracks in
the adhesive, which deposits some energy without causing a heat production
and which can be get back if healing is allowed. A specific example
might be}
\begin{subequations}\begin{align}
&a_0\big(
  \alpha,\thetaS
  \big):=
  -G_{_{\rm C}}^{}(\thetaS)\alpha
\intertext{with the so-called fracture energy (or {\it fracture toughness})
  $G_{_{\rm C}}^{}:
  \R^+\to\R^+$ continuous and bounded
\COL{ while the rate-dependent dissipation part leading to the heat production 
might be}}
&a_1\big(\REFF{\JUMP{u}{}},\COLOR{\alpha,}\thetaS,\DT\alpha\big)
:=\frac12\DT\alpha^2\begin{cases}
\epsilon_\text{\sc dam}^{}\big(\REFF{\JUMP{u}{}},\alpha,\thetaS\big)
&\text{if }\ \DT\alpha\le0\\
1/\epsilon_\text{\sc heal}^{}\big(\REFF{\JUMP{u}{}},\alpha,\thetaS\big)
&\text{if}\ \DT\alpha\ge0\end{cases}
\label{a1}
\end{align}\end{subequations}
\COL{with some small $\epsilon_\text{\sc dam}>0$ and $\epsilon_\text{\sc heal}>0$
  causing small heat production during fast delamination and very slow healing,
  respectively.
}
\end{remark}


\begin{remark}[\COL{Some variants of the model}.]\label{rem-others}
\upshape
Actually, there are several \COLOR{variations of the presented model}. 
E.g., considering the Signorini
contact \COLOR{(resulting by sending $\gammaC(\REFF{\JUMP{u}{\rmn}})\to+\infty$ 
if $\REFF{\JUMP{u}{\rmn}}<0$ in \eqref{normal-compliance})}
does not comply with \eqref{ass-gamma-c}
and needs either ignoring heat transfer or inertial effects, cf.\ e.g.\
\cite{RosRou11TARI,RosRou13ACDM,ScaSch17CPVB}.
\COL{Engineering models typically do not allow for healing, i.e.\
  $\epsilon_\text{\sc heal}=+\infty$ in \eqref{a1}, which needs} weakening of
the continuity assumption \eqref{ass:a1} \COL{but still}
allows for a limit passage in the conventional weak formulation
of the delamination flow rule, cf.\ again \cite{RosRou11TARI,RosRou13ACDM}.
If the adiabatic effects on $\GC$ would be suppressed (by setting $b_0=0$),
then 
the delamination $\alpha$ might be allowed to evolve rate-independently
by considering
\COL{
both $\epsilon_\text{\sc dam}=0$ and $\epsilon_\text{\sc heal}=+\infty$
in \eqref{a1}}.
\REFF{The viscosities in the adhesive and the bulk (cf.\ Remark~\ref{rem-RI-pi}
below) then avoid}
the question of a suitable concept of weak solutions.
\COL{Let us point out that, otherwise,} 
a combination of two rate-independent processes
governed not by a jointly convex potential \COL{leads to various quite distinct 
solution concepts}, cf.\ \cite{MieRou15RIS} for a thorough discussion;
actually our semi-implicit formula \eqref{system-disc} would lead to a certain
stress-driven-like weak solution, cf.\ \cite{vodicka14A1}.
%
\end{remark}


\REFF{
\begin{remark}[Hidden rate-dependency of plastic slip.]\label{rem-RI-pi} 
\upshape
The reason for considering viscosity in the adhesive (i.e.\ the
dumpers $d_{\rmn}$ and $d_{\rmt}$ in Fig.\,\ref{fig:rheology}) is to control 
the slip rate $\DT\pi$ in $L^2(\SC;\R^{d-1})$ even though there 
is no explicit rate dependence in the flow rule \eqref{adhes-form-d7};
note that $\JUMP{\DT u}{}$ is well controlled in $L^2(\SC;\R^{d-1})$
due to the viscosity in the bulk.
Otherwise,  $\DT\pi$ would be controlled only in $L^1(\SC;\R^{d-1})$ or, more
precisely, in measures on the closure of $\SC$, which would make
analytical troubles unless $\sigma_{\rm y}$ would not be constant or,
at least, unless $d=2$ and $\sigma_{\rm y}$ would depend only 
on $\alpha$, which is then in $C(\bar\SC)$, but not on $\thetaS$.
\end{remark}
}

\COL{
\begin{remark}[Thermally expansive adhesive.]
\upshape
Like the term $\thetaB\bbB:e$ in \eqref{psi-bulk},
the term $b_0(\alpha)\thetaS$ in \eqref{psi-surf} could be 
replaced by $\thetaS b_0(\alpha)\cdot\JUMP{u}{}$ with some 
$b_0:\R\to\R^{d-1}$. This modification could model a thermally expansive adhesive
while not substantially affecting the analysis.
\end{remark}


}

\REF{


\section{Extension towards poro-viscoelastic media and adhesives}\label{poro}

The set of internal variable (so far composed from the 
interfacial damage and plastic-like interfacial slip)
can be enhanced in particular applications. One possible way is
to consider damage and plastic strain also in the bulk $\varOmega_1\cup\varOmega_2$.
Another enhancement may consider also other physical processes both in
the bulk and in the adhesive on the interface. E.g.\ delamination of 
piezoelectric materials has been extensively studied in literature.

Here we want to illustrate the applicability of the structure 
\eqref{system-abstract} on processes with nonlocal dissipation potentials
and which are undergoing not only in the adhesive surface but possibly
also in the bulk. More specifically, both humidity/moisture 
(or water pressure) and 
temperature are know to affect the debonding processes in specific 
applications, cf.\ e.g.\ 
\cite{AlMaSa06CDFI,Asp98EMTI,BruSau95WFIC,KKMS17MASD,Shil17MDCV}
or also \cite[Chap.\,6]{Sih91MFIP}. Then,
beside heat transfer, also a water transport in the porous adhesive
layer and/or in the poro-viscoelastic bulk is to be considered. 
In fact, the diffusant need not be just water (as typically in 
civil engineering constructions or in rock or soil mechanics) but 
some solvents in polymers or hydrogen in metals undergoing the
so-called metal/hydride transition, cf.\ e.g.\ \cite{RouTom15TDMD}. 
In some situations, only solvent 
transport over the adhesive interface is relevant, while in other
situations the transport inside the bulk accompanied with pronounced swelling 
effects may lead to delamination, cf.\ Remark~\ref{rem-c=0} 
below. Of course, both mentioned situations
can occur simultaneously, which will be covered by the model presented below.

Models used in engineering are often very phenomenological and without
rational-mechanical structure. In contrast, consistently with the 
previous part, here we want to cast a thermodynamically consistent
model where the diffusion is governed by the chemical potential 
(containing also a pressure) so that we cover both Fick and Darcy 
law. Other concepts used are the Biot model for saturated flow \cite{Biot41GTTS}
and capillarity due to 
the gradient of concentration involved in the free energy, leading 
to the Cahn-Hilliard model \cite{CahHil58FEUS}.

We denoting by $\zetaS$ and $\zetaB$ the content of diffusant in the
interfacial adhesive layer and in the bulk, respectively. Then the specific
bulk free energy \eqref{psi-bulk} expands as:
\begin{align}\nonumber
\psi\indexV (e,\zetaB,\thetaB,\nabla\zetaB)&=
\frac12\bbC e{:} e-(\thetaB{-}\theta_{_{\rm R}})\bbB{:}e-{\psi\indexV}_0(\thetaB)
\\[-.3em]
\label{psi-bulk+}
&\ \ \ \ +\frac12\MB(\betaB{\rm tr}\,e-\zetaB)^2+\frac12K(\zetaB{-}
\zetaBeq)^2
+\frac\kappa2|\nabla\zetaB|^2
\end{align}
with $\bbB=\bbC\bbE$ as in \eqref{psi-bulk} and with $M>0$ the so-called 
Biot modulus, $\betaB>0$ the Biot coefficient, $K>0$ a coefficient (related
but not identical with the bulk modulus) and ${\zetaBeq}>0$ the 
equilibrium content.

Analogously we also cast the interface model. In contrast to \eqref{psi-surf},
the analog of the Biot term couples the normal and the tangential components
of $\JUMP{u}{}$ because the analog of the strain trace, i.e.\ 
${\rm tr}\,e={\mathbb I}{:}e$ with ${\mathbb I}\in\R^{d\times d}$ denoting the 
identity matrix, is the sum of components of $\JUMP{u}{}$, i.e.\ 
$\JUMP{u}{\rmn}+{\rm I}_{d-1}\cdot
\JUMP{u}{\rmt}={\rm I}_d\cdot\JUMP{u}{}$ with the vector 
${\rm I}_d=(1,...,1)\in\R^d$. Then the specific contact interface energy
\eqref{psi-surf} is enhanced as
\begin{align}\nonumber
&\psi\indexS\big(\JUMP{u}{},
\pi,\alpha,\zetaS,\thetaS,\nablaS\pi,\nablaS\alpha,\nablaS\zetaS\big)=
\int_{\GC}
\Big(\gamma_{\rmt}\big(\alpha,\JUMP{u}{\rmt}\!{-}\pi\big)
+\gamma_{\rmn}\big(\alpha,\JUMP{u}{\rmn}\big)+\gammaC\big(\JUMP{u}{\rmn}\big)
\\&\quad\nonumber
+\frac12\MS\big(\betaS\big({\rm I}_{d-1}\cdot(\JUMP{u}{\rmt}\!{-}\pi)
+\JUMP{u}{\rmn}\big)-\zetaS\big)^2
+\frac12\KS(\zetaS{-}{\zetaSeq})^2
+\frac{\kappa_{\rm H}}2|\pi|^2+\frac{\kappa_1}2|\nablaS\pi|^2
\\&\quad
+\frac{\kappa_2}2|\nablaS\alpha|^2+\frac{\kappa_3}2|\nablaS\zetaS|^2
-b_0(\alpha)\thetaS-{\psi\indexS}_{\,0}(\thetaS)\Big)\,\d S
+\begin{cases}-\,a_0(\alpha)&\text{if $\ 0\!\le\!\alpha\!\le\!1$ a.e.~on $\GC$},
\\\ \,+\infty&\text{otherwise}.\end{cases}
\label{psi-surf+}\end{align}
Like in \eqref{Phi-epos}, the overall free energy is 
\begin{align}\label{Phi-epos+}
\varPsi(e,\JUMP{u}{},\pi,\alpha,\thetaS,\thetaB)=
\int_{\GC}\!\!\psi\indexS\big(\JUMP{u}{},
\pi,\alpha,\zetaS,\thetaS,\nablaS\pi,\nablaS\alpha,\nablaS\zetaS\big)\,\d S
+\int_{\varOmega {\setminus} \GC}\!\!
\psi\indexV(e,\zetaB,\thetaB,\nabla\zetaB)\, \d x\,.
\end{align}
The driving force for diffusion is the gradient of the chemical potential which
is a derivative of the free energy according concentration. Here,
we have both the bulk and the adhesive chemical energies, denoted
respectively by $\muB$ and $\muS$. Because of the gradients of $\zetaB$ and of $\zetaS$,
we must consider the functional derivative, which gives
\begin{subequations}\label{chem-pot}\begin{align}
\muS&=\partial_{\zetaS}\!\varPsi=\MS(\zetaS{-}\betaS{\rm I}_d{\cdot}\JUMP{u}{})
+K(\zetaS{-}{\zetaSeq})-{\divS}(\kappa_3\nablaS\zetaS)\,,
\\\muB\,&=\partial_c\varPsi\,=\MB(\zetaB-\betaB{\rm tr}\,e(u))+K(\zetaB-{\zetaBeq})-{\rm div}(\kappa\nabla\zetaB)\,.
\end{align}\end{subequations}
The diffusion equations in the bulk and in the adhesive layer then
represent conservation of the mass of the diffusant:
\begin{subequations}\label{diffusion}\begin{align}
&\label{diffusion2}
\DTzetaS-\divS(\bbMS(\alpha,\zetaS,\thetaS)\nablaS\muS)
=m(\muB|_{\varOmega_1}^{}+\muB|_{\varOmega_2}^{}-2\muS)\,,
\\\label{diffusion1}
&\DT\zetaB-{\rm div}(\bbMB(\zetaB,\thetaB)\nabla\muB)=0\,,
\intertext{where $\bbMS:\R^3\to\R^{(d-1)\times(d-1)}$ and 
$\bbMB:\R^2\to\R^{d\times d}$ are the mobility 
tensors in the adhesive and the bulk, respectively, 
and $m=m(\alpha)$ is a mobility coefficient between the bulk and the adhesive.
 Of course, it should be accompanied by some boundary conditions,
for simplicity considered for an isolated body:}
&\label{diffusion3}
(\bbMB(\zetaB,\thetaB)\nabla\muB|_{\varOmega_i}^{})\cdot\vec{n}=\begin{cases}0&\text{ on 
$(I\times\partial\varOmega_i)\setminus\SC$},
\\m(\muS{-}\muB|_{\varOmega_i}^{})&\text{ on $\SC$},\ \ i=1,2,\end{cases},
\\&
\nablaS\muS\cdot\nu\indexS=0\hspace{4em}\text{on $\partial\SC$}.
\end{align}\end{subequations}
The right-hand side of \eqref{diffusion2} is the source/sink of the
diffusant in the adhesive layer due to the transfer from the adjacent
bulk region, and $\muB|_{\varOmega_1}$ means the trace of $\muB$ in
$I\times\varOmega_1$ on $\SC$ and analogously for $\muB|_{\varOmega_2}$.
This flux occurs naturally in the boundary condition
\eqref{diffusion3} on $\SC$. In principle, the coefficient $m$ may
depend on $\alpha$, $\thetaS$, and on the traces of $\thetaB|_{\varOmega_1}$ 
$\thetaB|_{\varOmega_2}$ on $\SC$, or may even be different on the boundary of 
$\varOmega_1$ and $\varOmega_2$. Of course,
this system is to be completed by the initial conditions, namely
\begin{align}\label{diffusion-IC}
\zetaS(0)={\zetaS}_{,0}\quad\text{and}\quad\zetaB(0)={\zetaB}_{,0}\,.
\end{align}

This diffusion problem \eqref{diffusion} is coupled with the original model
thermomechanical model \eqref{bulk-model} and its boundary 
 condition \eqref{BC} as well as the contact interfacial model \eqref{system} 
by considering the stress $\sigma=(\psi\indexV)_e'$
taking into account \eqref{psi-bulk+}, which now gives 
\begin{align}\label{stress+}
\sigma=\bbD
e(\DT{u})+\bbC(e(u){-}\bbE\thetaB)+
\frac1d{\betaB}^2\MB\Big({\rm sph}\,e(u){-}\frac1d\zetaB{\mathbb I}\Big)
\end{align}
with ${\rm sph}\,e=({\rm tr}\,e){\mathbb I}/d$ denoting the spherical (sometimes
called volumetric) part of the strain tensor $e$. This is now to be used
in \eqref{adhes-form-d1}, and thus in \eqref{eq6:adhes-class-form1}
and in \eqref{eq6:adhes-class-form3-bis}, too.


Consistently with the previous 
notation, we denote $\bfmu=(\muS,\muB)$ and $\bfzeta=(\zetaS,\zetaB)$
and formulate the diffusion problem \eqref{diffusion} in a monolithic way.
Again we use the gradient/difference operator $\mathcal{G}$ from  
\eqref{def-of-G} and, like $\mathcal{K}(\bftheta)$, we define the monolithic
mobility operator 
\begin{align}\label{def-of-M}
&
\mathcal{M}(\alpha,\bfzeta):=\big(\bbM\indexS(\alpha),
m(\alpha,\thetaS),
\bbM\indexV(\zetaB)\big)
\end{align}
Then the diffusion problem \eqref{chem-pot}--\eqref{diffusion} can be written 
simply as
\begin{align}\label{diffusion-compact}
\DT\bfzeta+\mathcal{G}^*\big(\mathcal{M}(\alpha,\bfzeta)\mathcal{G}\bfmu\big)
=0\ \ \ \text{ with }\ \ \
\bfmu=\partial_{\bfzeta}\varPsi(\boldsymbol{E}u,\pi,\bfzeta)\,,
\end{align}
where we used already the special choice \eqref{psi-surf+}
which causes, in particular, that $\partial_{\bfzeta}\varPsi$ does not
depend on $\alpha$ and $\bftheta$.

The energetics of this model
\eqref{diffusion-compact} coupled with (\ref{system-abstract})
can be revealed by testing the equations in
\eqref{diffusion-compact} respectively by $\bfmu$ and $\DT\bfzeta$,
and the equations in (\ref{system-abstract}a,b,c,d) again respectively by
$\DT u$, $\DT\pi$, $\DT\alpha$, and $\one$.
%
Instead of \eqref{Psi}, the overall free energy $\varPsi$ now is
from \eqref{Phi-epos+}, while the overall dissipation rate \eqref{dissip-rate}
is enhanced by an extra-dissipation due to diffusion processes so that
altogether it reads as:
\begin{align}\nonumber
&R\big(\JUMP{u}{},\alpha,\zetaS,\zetaB,\thetaS,\thetaB;\DT{e},
  \JUMP{\DT u}{},\DT\pi,\DT{\alpha},\nablaS\muS,\muS{-}\muB|_{\varOmega_i}^{},
  \nabla\muB\big)
\\&\nonumber\qquad:=\int_{\GC}\!\xi\indexS\big(\JUMP{u}{},\alpha,\thetaS;
\JUMP{\DT u}{},\DT\pi,\DT{\alpha}\big)
+\bbMS(\zetaS,\thetaS)\nablaS\muS{\cdot}\nablaS\muS
+\sum_{i=1}^2m\big(\muS{-}\muB|_{\varOmega_i}^{}\big)^2\d S
\\[-.4em]&\qquad\qquad\qquad\qquad\label{dissip-rate+}
+\int_{\varOmega{\setminus}\GC}\!\!\xi\indexV(\thetaB;\DT{e})+
\bbMB(\zetaB,\thetaB)\nabla\muB{\cdot}\nabla\muB\,\d x\,.
\end{align}
Written monolithically, it means
\begin{align}\nonumber
R\big(\JUMP{u}{},\alpha,\bfzeta,\bftheta;
\boldsymbol{E}\DT u,\DT\pi,\DT{\alpha},\mathcal{G}\bfmu\big)
=&\int_{\GC}\!\xi\indexS\big(\JUMP{u}{},\alpha,\thetaS;
\JUMP{\DT u}{},\DT\pi,\DT{\alpha}\big)\,\d S
\\[-.4em]&\ 
+\int_{\varOmega{\setminus}\GC}\!\!\xi\indexV(\thetaB;\DT{e})\,\d x
+
\langle\hspace*{-.25em}\langle\mathcal{M}(\alpha,\bfzeta)\mathcal{G}\bfmu\,\,
\mathcal{G}\bfmu\rangle\hspace*{-.25em}\rangle\,,
\end{align}
where the bilinear form
$\langle\hspace*{-.25em}\langle\cdot,\cdot\rangle\hspace*{-.25em}\rangle$ is
again from \eqref{<<.>>}. The specific dissipation rates $r\indexS$ and
$r\indexV$ in (\ref{entropy-eq}a,b) are then expanded by the interfacial heat
source $\bbMS(\zetaS,\thetaS)\nablaS\muS{\cdot}\nablaS\muS$ and by
the bulk heat source $\bbMB(\zetaB,\thetaB)\nabla\muB{\cdot}\nabla\muB$,
and also part (say 1/2) of the heat production due to the last term in 
\eqref{dissip-rate+}, so that 
\begin{subequations}\label{heat-enhanced}\begin{align}
&\nonumber
r\indexS=
\friction(\alpha,\thetaS)\gammaCprime\big({\JUMP{u}{}}_\rmn\big)
\big|\JUMP{\DT u}{\rmt}\big|
+\sigma_{\rm y}(\alpha,\thetaS)|\DT\pi|
+\DT\alpha\partial_{\DT\alpha}a_1\big(\JUMP{u}{},\alpha,\thetaS,\DT\alpha\big)
+d_\rmt(\alpha,\thetaS)\big|\JUMP{\DT u}{\rmt}\!{-}\DT\pi\big|^2
\\[-.3em]&\qquad\qquad
+d_\rmn(\alpha,\thetaS)\JUMP{\DT u}{\rmn}^{\!\!\!\!2}+\bbMS(\zetaS,\thetaS)
\nablaS\muS{\cdot}\nablaS\muS
+\frac12\sum_{i=1}^2m\big(\muS{-}\muB|_{\varOmega_i}^{}\big)^2,\ \ \text{ and}
\\
&r\indexV=\bbD e(\DT{u}){:}e(\DT{u})
+\bbMB(\zetaB,\thetaB)\nabla\muB{\cdot}\nabla\muB,
\end{align}\end{subequations}
while the resting heat production due to the last term in 
\eqref{dissip-rate+} should go into the adjacent bulk domains through 
enhancement of the right-hand side in the boundary condition 
\eqref{BC-heat-GC} by $\frac12(\muS{-}\muB|_{\varOmega_i}^{})^2$.

Actually, one expects the dissipation rate to depend on rates,
i.e.\ here on $\DT\zetaS$ and $\DT\zetaB$ rather than on the gradients 
$\nablaS\muS$ and 
$\nabla\muB$ or the difference $\muS{-}\muB|_{\varOmega_i}^{}$.
Indeed, the chemo-mechanical part of the problem essentially enjoys 
the general structure \eqref{system-abstract} with the rate-dependent 
(pseudo)potential of dissipative forces which is now nonlocal, however. 
To see this quite well-known fact, we adapt a bit technical calculation 
\cite{Roub17ECTD} to our monolithic bulk/interfacial model. We realize the 
structure of \eqref{diffusion-compact}
\begin{align}\label{diffusion-compact+}
  \DT\bfzeta+[\varXi_\mathcal{M}^{\rm dif}]'(\bfmu)=0\ \ \ \text{ with }
\ \ \ \bfmu=\partial_{\bfzeta}\varPsi(\boldsymbol{E}u,\pi,\bfzeta)\,,
\end{align}
with $\mathcal{M}=\mathcal{M}(\alpha,\bfzeta)$ from \eqref{def-of-M}, and
where $\varXi_\mathcal{M}^{\rm dif}$ denotes a quadratic potential
\begin{align}\nonumber
\varXi_\mathcal{M}^{\rm dif}:\bfmu=(\muS,\muB)\mapsto\frac12
\int_{\GC}\bbMS\nablaS\muS{\cdot}\nablaS\muS+
\frac12\sum_{i=1}^2m\big(\muS{-}\muB|_{\varOmega_i}^{}\big)^2\,\d S\qquad\qquad
\\[-.4em]+
\frac12\int_{\varOmega\setminus\GC}\!\!\!
\bbMB\nabla\muB{\cdot}\nabla\muB\,\d x
=\frac12\langle\hspace*{-.25em}\langle\mathcal{M}\mathcal{G}\bfmu\,\,
\mathcal{G}\bfmu\rangle\hspace*{-.25em}\rangle\,.
\label{diffusion-pot}\end{align}
Thus, the G\^ateaux derivative of $\varXi_\mathcal{M}^{\rm dif}$ in 
\eqref{diffusion-compact+}
``monolithically'' expresses the combination of the Beltrami-Laplace operator
with the mobility matrix $\bbMS$  on $\GC$ and the Laplace operator with the
mobility matrix $\bbMB$ on $\varOmega\setminus\GC$, and the 
transition conditions between $\GC$ and the adjacent bulk. More specifically,
$\Delta_\mathcal{M}$ is the linear positive-definite operator 
on $H^1_\equiv$ denoting 
the Hilbert space $H^1(\GC)\times H^1(\varOmega\setminus\GC)$ modulo
constant functions, or essentially equivalently 
$H^1_\equiv:=\{\bfmu=(\muS,\muB)\in H^1(\GC)\times
H^1(\varOmega{\setminus}\GC);\ \int_{\GC}\muS\,\d x+\int_\varOmega\muB\,\d x=0\}$. 
Here $\bbMS$ and $\bbMB$ are assumed symmetric to ensure existence of
a potential of this linear operator.

The convex conjugate $[\varXi_\mathcal{M}^{\rm dif}]^*$ of
$\varXi_\mathcal{M}^{\rm dif}$ is also quadratic 
and defines by means of its G\^ateaux differential the inverse
operator $[\varXi_\mathcal{M}^{\rm dif}]'$ which is nonlocal. 
Realizing the definition of the chemical potentials, \eqref{diffusion-compact+} 
then rewrites as 
\begin{align}\label{diffusion-compact++}
\big[[\varXi_\mathcal{M}^{\rm dif}]^{*}\big]'\DT{\bfzeta}
=\partial_{\bfzeta
}\varPsi(\boldsymbol{E}u,\pi,\bfzeta)\,,
\end{align}
which already fits with the form of \eqref{system-abstract}.
The dissipation rate \eqref{dissip-rate+} can be then equivalently 
expressed in terms of rates as 
\begin{align}\nonumber
  &\!R\big(\JUMP{u}{},\alpha,\bfzeta,\bftheta;
  \boldsymbol{E}\DT u\DT\pi,\DT{\alpha},\DT\bfzeta\big)
:=
\int_{\GC}\!\xi\indexS\big(\JUMP{u}{},\alpha,\thetaS;
\JUMP{\DT u}{},\DT\pi,\DT{\alpha}\big)\d S
\\[-.4em]&\hspace{14em}
\label{dissip-rate++}
+\!\int_{\varOmega{\setminus}\GC}\!\!\!\!\xi\indexV(\thetaB;\DT{e})\,\d x
+2\big[\varXi_{
\mathcal{M}(\alpha,\bfzeta)}^{\rm dif}\big]^*(\DT\bfzeta)\,,\!
\end{align}
where $[\cdot]^*$ denotes the mentioned convex conjugate functional.
 
The entropy balance \eqref{entropy-overall} now uses the enhanced 
heat production \eqref{heat-enhanced} together with the entropy
production $\sum_{i=1,2}\int_{\varOmega_i}\frac12(\muS{-}\muB|_{\varOmega_i}^{})^2/{\thetaB}_i$
on the boundaries $\partial\varOmega_i\cup\GC$ due to the mentioned additional 
heat production in the boundary conditions \eqref{BC-heat-GC}.

As the free energy is now convex (even quadratic) in
$(u,\JUMP{u}{},\pi,\zetaS,\zetaB)$, i.e.\ (\ref{sep-convex}a,b) now reads as
\begin{subequations}\label{sep-convex+}\begin{align}\label{sep-convex1+}
&(u,\pi,\bfzeta)\mapsto\varPsi(\boldsymbol{E}u,\pi,\alpha,\bfzeta,\bftheta)\ \ \text{ is convex for all
$\alpha$ and $\bftheta$},
\\&\alpha\mapsto\varPsi(\boldsymbol{E}u,\pi,\alpha,\bfzeta,\bftheta)\ \ \ \text{ is convex for all
$u$, $\pi$, and $\bftheta$},
\end{align}\end{subequations}
the system \eqref{system-abstract} enhanced by the diffusion equation
\eqref{diffusion-compact+}
bears the energy-conserving discretisation \cite{Roub17ECTD}
similarly as \eqref{system-disc} as a three-step decoupled scheme.
The first sub-problem now handles
$(\boldsymbol{E}u_\tau^k,v_\tau^k,\pi_\tau^k,\bfzeta_\tau^k)$. 
More specifically, having in mind the monolithic form
\eqref{diffusion-compact++}, thanks to \eqref{sep-convex1+}, the first step
(\ref{system-disc}a-c) will be now made jointly with 
\begin{subequations}\label{diffusion-disc}\begin{align}
&\frac{\bfzeta_\tau^k{-}\bfzeta_\tau^{k-1}}\tau+
\big[\varXi_{\mathcal{M}(\alpha_\tau^{k-1},\bfzeta_\tau^{k-1})}\big]'\bfmu_\tau^{k-1/2}=0,
\label{diffusion-disc1}\\&
\bfmu_\tau^{k-1/2}:=\frac{\bfmu_\tau^{k}+\bfmu_\tau^{k-1}}2=
\partial_{\bfzeta}\varPsi(\boldsymbol{E}u_\tau^{k-1/2},\pi_\tau^{k-1/2},\bfzeta_\tau^{k-1/2})\,,
\label{diffusion-disc2}
\end{align}\end{subequations}
where that the difference formula 
\eqref{def-of-circ} simplifies to  \eqref{diffusion-disc2}
because $\varPsi$ is quadratic in the variable involved in 
$\partial_{\bfzeta}\varPsi$. Of course, the recursive scheme 
\eqref{diffusion-disc} should start for $k=1$ by using the initial conditions 
$\bfzeta_0=({\zetaS}_{,0},{\zetaB}_{,0})$.
Like \eqref{potential-for-upi}, the resulted recursive boundary-value problems
have potentials. More specifically, $\varPsi$ in 
\eqref{potential-for-upi} now is to contain also the argument 
$(\bfzeta{+}\bfzeta_\tau^{k-1})/2$ and the whole functional
\eqref{potential-for-upi} is to be expanded by the term
$[\varXi_{\mathcal{M}(\alpha,\bfzeta)}^{\rm dif}]^*((\bfzeta{-}\bfzeta_\tau^{k-1})/\tau)$,
cf.\ also \cite[Formula (17a)]{Roub17ECTD} for the bulk part.

The energetics and the a-priori estimates  related to this diffusion
enhancement can be obtained by testing \eqref{diffusion-disc1} by
$\bfmu_\tau^{k-1/2}$ and \eqref{diffusion-disc2} by
$(\bfzeta_\tau^k{-}\bfzeta_\tau^{k-1})/\tau$, which yields

\begin{align}\nonumber
  0&=\big\langle
  \big[\varXi_{\mathcal{M}(\alpha_\tau^{k-1},\bfzeta_\tau^{k-1})}\big]'\bfmu_\tau^{k-1/2},
  \bfmu_\tau^{k-1/2}\big\rangle+
  \Big\langle\partial_{\bfzeta}\varPsi(\boldsymbol{E}u_\tau^{k-1/2},\pi_\tau^{k-1/2},\bfzeta_\tau^{k-1/2}),\frac{\bfzeta_\tau^k{-}\bfzeta_\tau^{k-1}}\tau\Big\rangle
  \\&\nonumber
  =2\varXi_{\mathcal{M}(\alpha_\tau^{k-1},\bfzeta_\tau^{k-1})}(\bfmu_\tau^{k-1/2})
  +\frac{\varPsi(\boldsymbol{E}u_\tau^{k-1/2},\pi_\tau^{k-1/2},\bfzeta_\tau^{k})
  -\varPsi(\boldsymbol{E}u_\tau^{k-1/2},\pi_\tau^{k-1/2},\bfzeta_\tau^{k-1})}\tau\,,
\end{align}
which is to be added to \eqref{system-disc1-3-tested} and then to
\eqref{system-disc1-4-tested}, too.
In addition to the a-priori estimates \eqref{est}, assuming $\bbMS$ and 
$\bbMS$ positive definite, $\MS>0$, $\MB>0$, and $m\ge0$, and also 
${\zetaS}_{,0}\in H^1(\GC)$ and ${\zetaB}_{,0}\in H^1(\varOmega\setminus\GC)$, we 
now obtain also  
\begin{subequations}\label{est+}\begin{align}
&
\big\|\nablaS{\muS}_\tau\big\|_{L^2(\SC;\R^{d-1})}\le C
\ \ \ \text{ and }\ \ \ 
\big\|\nabla{\muB}_{,\tau}\big\|_{L^2(Q;\varOmega\setminus\SC;\R^d)}\le C\,,
\\
&\big\|{\zetaS}_{,\tau}\big\|_{L^\infty(0,T;H^1(\GC))\,\cap\,
H^1(0,T;H^1(\GC)^*)}^{}\le C\ \ \text{ and }\ \ 
\big\|\divS(\kappa_3\nablaS{\zetaS}_{,\tau})\big\|_{L^\infty(0,T;L^2(\GC))}\le C\,,
\\
&\big\|{\zetaB}_{,\tau}\big\|_{L^\infty(0,T;H^1(\varOmega\setminus\GC))\,\cap\,
H^1(0,T;H^1(\varOmega\setminus\GC)^*)}\le C
\ \ \ \text{ and }\ \ \ 
\big\|{\rm div}(\kappa\nabla{\zetaB}_{,\tau})\big\|_{L^\infty(0,T;L^2(\varOmega))}^{}\le C\,.
\end{align}\end{subequations}
This guarantees numerical stability of the staggered discretisation
enhanced by diffusion, and also allows for convergence towards
the weak solution now enhanced by suitable weak formulation of 
the initial-boundary-value problem 
\eqref{chem-pot}--\eqref{diffusion}--\eqref{diffusion-IC},
cf.\ \cite{Roub17ECTD} for details as far as the bulk model concerns
while the $(d{-}1)$-interfacial variant is analogous.

\begin{remark}[More coupling.]\upshape
The dissipation mechanism in the thermo-mechanical part can easily be
made $(\zetaS,\zetaB)$ dependent. For example, the friction coefficient 
$\friction$ often depend on the diffusant content $\zetaS$. Also, the 
coefficient $m$ can depend on $\alpha$ and on $\JUMP{u}{}$. In some 
applications it is also reasonable to make elastic moduli in the bulk 
$\bbC$ or in the adhesive $\gamma_\rmt(\alpha,\cdot)$ and 
$\gamma_\rmn(\alpha,\cdot)$ dependent on the diffusant content $\zetaS$ or 
$\zetaB$, respectively. Yet, this latter coupling would affect the 
chemical potentials $\muS$ and $\muB$ and might make the analysis and 
computer algorithmic implementation more complicated.
\end{remark}

\begin{remark}[Models with $\zetaB=0$.]\label{rem-c=0}\upshape
  Some engineering models consider moisture propagation only along $\GC$
  in the adhesive, while the bulk is pores-free. This situation is
  covered by the above model simply by putting the initial bulk content
  ${\zetaB}_{,0}=0$ and the transient adhesive-bulk coefficient $m=0$.
\end{remark}


}


\subsubsection*{Acknowledgment}
\COL{The author is deeply thankful for many inspiring discussions with 
Vladislav Manti\v{c} and Roman Vodi\v cka.
Also, many comments and conceptual suggestions of two anonymous referees
are truly acknowledged. This research has been partly covered} by the Czech
Science Foundation through the grants 
16-34894L ``Variational 
structures in continuum thermomechanics of solids'' and
17-04301S ``Advanced mathematical methods for dissipative evolutionary 
systems''.
Eventually, the hospitality of the University 
of Seville \COL{within several} stays during  2015--\COL{2018
as well as} the institutional support RVO: 61388998 (\v CR)
is acknowledged, too.

\end{document}